\documentclass{amsart}

\usepackage{color}
\usepackage{soul}
\usepackage{amsmath, amssymb}

\usepackage{hyperref}
\hypersetup{
    colorlinks=true, 
    linktoc=all,     
    linkcolor=blue  
}

\newcommand{\dotieconcat}[2]{
  \text{\raisebox{.8ex}{$\smallfrown$}}%
}

\newcommand{\C}{\mathbb{C}}

\newcommand{\N}{\mathbb{N}}
\newcommand{\F}{\mathbb{F}}

\newcommand{\norm}[1]{\left\| #1 \right\|}

\newcommand{\boldB}{\mathbf{B}}

\newcommand{\boldD}{\mathbf{D}}

\newcommand{\cstar}{$\mathrm{C}^*$}

\newbox\gnBoxA
\newdimen\gnCornerHgt
\setbox\gnBoxA=\hbox{$\ulcorner$}
\global\gnCornerHgt=\ht\gnBoxA
\newdimen\gnArgHgt
\def\Quinequote #1{%
    \setbox\gnBoxA=\hbox{$#1$}%
    \gnArgHgt=\ht\gnBoxA%
    \ifnum     \gnArgHgt<\gnCornerHgt \gnArgHgt=0pt%
    \else \advance \gnArgHgt by -\gnCornerHgt%
    \fi \raise\gnArgHgt\hbox{$\ulcorner$} \box\gnBoxA %
    \raise\gnArgHgt\hbox{$\urcorner$}}

\DeclareMathOperator{\Tr}{Tr}
\DeclareMathOperator{\tr}{tr}
\DeclareMathOperator{\Rea}{Re}
\DeclareMathOperator{\spec}{sp}
\newcommand{\UCP}{\mathrm{UCP}}

\newcommand{\lmax}{\lambda_{\max}}
\newcommand{\Mn}{M_n(\mathbb{C})}
\newcommand{\va}{\vec{a}}
\newcommand{\vb}{\vec{b}}
\newcommand{\vc}{\vec{c}}
\newcommand{\vx}{\vec{x}}
\newcommand{\HS}{\mathrm{HS}}

\theoremstyle{plain}
\newtheorem{theorem}{Theorem}[section]
\newtheorem{lemma}[theorem]{Lemma}
\newtheorem{corollary}[theorem]{Corollary}
\newtheorem{proposition}[theorem]{Proposition}

\newtheorem{question}[theorem]{Question}
\newtheorem{fact}[theorem]{Fact}

\newtheorem*{theoremA}{Theorem A}
\newtheorem*{theoremB}{Theorem B}
\newtheorem*{theoremC}{Theorem C}
\newtheorem*{theoremD}{Theorem D}
\newtheorem*{theoremE}{Theorem E}
\newtheorem*{theoremF}{Theorem F}
\newtheorem*{theoremG}{Theorem G}
\theoremstyle{definition}
\newtheorem*{questionA}{Question A}

\theoremstyle{definition}
\newtheorem{definition}[theorem]{Definition}

\newtheorem{examples}[theorem]{Examples}
\newtheorem{remark}[theorem]{Remark}

\numberwithin{equation}{section}

\begin{document}
\title{Effectively nuclear \cstar-algebras}
\author[C. J. Eagle]{Christopher J. Eagle${}^1$} 
\address[C. J. Eagle]
{Department of Mathematics and Statistics, University of Victoria. PO BOX 1700 STN CSC, Victoria, British Columbia, Canada. V8W 2Y2}%
\email{eaglec@uvic.ca}
\urladdr{http://www.math.uvic.ca/~eaglec}
\thanks{${}^1$ Supported by NSERC Discovery Grant RGPIN-2021-02459}

\author[I. Goldbring]{Isaac Goldbring${}^2$}
\address[I. Goldbring]{University of California, Irvine}
\email{igoldbri@uci.edu}
\urladdr{https://www.math.uci.edu/~isaac}
\thanks{${}^2$  Partially supported by NSF grant DMS-2054477.}

\author[T. H. McNicholl]{Timothy H. McNicholl}
\address[T. H. McNicholl]{Iowa State University}
\email{mcnichol@iastate.edu}
\urladdr{https://faculty.sites.iastate.edu/mcnichol/}

\author[R. Miller]{Russell Miller${}^3$}
\address[R. Miller]{Queens College, 65-30 Kissena Blvd, Queens, NY 11367 \& CUNY Graduate Center, 365 Fifth Avenue, New York, NY 10016, U.S.A.}
\email{russell.miller@qc.cuny.edu}
\thanks{${}^3$  Supported by NSF grant DMS-2348891,
by Simons Foundation award \#MP-TSM- 00007933,
and by several grants from the City University of New York PSC-CUNY Research Award Program.}

\begin{abstract}
We study the computable content of nuclearity for $C^*$-algebras. We show that
nuclearity effectivizes automatically: a presentation of a $C^*$-algebra is
effectively nuclear, in the sense that the completely positive approximations
witnessing nuclearity can be computed uniformly from the data they approximate,
if and only if it is a computable presentation of a nuclear $C^*$-algebra, and
likewise for the effectivization of the Choi--Effros form of nuclearity. The main
tool, which we believe to be of independent interest, is an effective version of
the Arveson extension theorem: we prove that the $n$-th matrix numerical range of
a tuple from a computable presentation is computably compact, uniformly in all of
its data, and deduce an algorithm which converts an approximately attainable tuple
of matrices into a computable unital completely positive map attaining it. A
byproduct is that the computable states of any computable presentation of a unital
$C^*$-algebra are weak$^*$-dense in its state space.

We also prove that the Cuntz algebras $\mathcal{O}_n$ are computably categorical
and obtain an effective form of Kirchberg--Phillips absorption. By contrast with
the conditions above, we do not know whether every computable presentation of a
nuclear $C^*$-algebra admits a computable embedding into $\mathcal{O}_2$ with a
computable conditional expectation; we answer this affirmatively for commutative
$C^*$-algebras, for AF algebras, and, by a dynamical argument using invariant
states, for the group $C^*$-algebras of a class of amenable groups including all
finitely presented residually finite ones. Finally, we show that the set of indices
of computable presentations of nuclear $C^*$-algebras is $\Pi^0_2$, and that the
corresponding set for c.e. presentations is $\Pi^0_2$-hard; the latter rests on a
result of independent interest, namely that amenability, and also inner
amenability, of computable groups is $\Pi^0_2$-complete within the class of
computable groups.
\end{abstract}
\maketitle


\section{Introduction}\label{sec:Intro}
 
A \cstar-algebra $A$ is \emph{nuclear} if the identity map on $A$ is the point-norm limit of
maps that factor through finite-dimensional matrix algebras by unital completely positive (ucp)
maps, that is, for every finite tuple $\vec a$ from $A$ and every $\epsilon>0$, there are $n \ge 1$
and ucp maps $\phi : A \to M_n(\mathbb{C})$ and $\psi : M_n(\mathbb{C}) \to A$ satisfying
$\max_i \|(\psi\circ\phi)(a_i) - a_i\| < \epsilon$. Nuclearity is one of the central finiteness
conditions in operator algebras: it is equivalent to the uniqueness of the \cstar-norm on algebraic
tensor products $A \odot B$, it is the \cstar-algebraic analogue of amenability for groups (a discrete
group $G$ being amenable exactly when $C^*_r(G)$ is nuclear), and it is a standing hypothesis in the
Elliott classification programme.
 
This paper is concerned with the computable content of nuclearity. The setting is the computable
structure theory of \cstar-algebras developed by Fox \cite{fox2022computable}, in which a \emph{presentation} $A^\#$ of a
\cstar-algebra $A$ consists of a distinguished sequence of points whose rational $*$-polynomial
combinations are dense; the presentation is \emph{c.e.} if the norms of these generated points are
uniformly right-c.e., and \emph{computable} if they are uniformly computable. In our previous
papers \cite{AFPaper,UHFPaper}, we began a programme of computable $K$-theory for \cstar-algebras, carried out in detail
for UHF and AF algebras; several of the tools developed there-the uniformly computable
presentations of matrix amplifications \cite[Proposition 1.7]{UHFPaper}, the passage from computably weakly
stable relations to c.e.\ closed sets \cite[Theorem 1.12]{UHFPaper}, and computable AF-certificates \cite[Main Theorem 1]{AFPaper}-are used throughout the present paper. The index set computations of \cite{AFPaper},
showing that the set of indices of computable presentations of UHF algebras, and likewise of (simple) AF
algebras, are $\Pi^0_2$-complete, are the direct precursors of Section \ref{sec:completeness} below.
 
The basic question we address is: \emph{what does it mean for a presentation to be nuclear in an
effective way, and is this stronger than being a computable presentation of a nuclear algebra?}  Given that our earlier works considered very special subclasses of the class of nuclear algebras, zooming out and considering the broader class as a whole is a natural next step.
There are several reasonable ways to effectivize the definition above. The most direct is to ask
that the ucp maps witnessing the completely positive approximation property be computable and be
produced by an algorithm; we call such a presentation \emph{effectively nuclear} (Definition \ref{def effnuc}(1)).
One may instead effectivize the Choi--Effros form of nuclearity, asking for a uniformly computable
sequence of finite-rank ucp maps $A^\# \to A^\#$ with computably bounded ranks converging computably
pointwise to the identity; we call this \emph{effectively Choi--Effros nuclear} (Definition \ref{def effnuc}(2)).
Our first main result is that none of these distinctions is real: nuclearity effectivizes
automatically.
 
\begin{theoremA}[Theorem \ref{mainequivalence}]
Let $A^\#$ be a presentation of a \cstar-algebra $A$. The following are equivalent:
\begin{enumerate}
\item $A^\#$ is computable and $A$ is nuclear;
\item $A^\#$ is effectively nuclear;
\item $A^\#$ is computable and effectively Choi--Effros nuclear.
\end{enumerate}
Moreover, the passage from $(1)$ to $(2)$ is uniform in a unital index for $A^\#$.
\end{theoremA}
 
We note that
the implication $(2) \Rightarrow (1)$ contains the observation that
effective nuclearity already forces the presentation to be computable (Lemma \ref{lem:automaticallycomputable}). It is natural to
compare Theorem A with the situation for groups, where the analogous phenomenon holds: a computable
presentation of a group admits a computable F{\o}lner sequence precisely when the group is amenable (see \cite{DudaIvanov}).
 
The engine behind Theorem A is an effective form of the Arveson extension theorem, which we
believe is of independent interest and which occupies Section \ref{sec:arveson}. Given a tuple
$\vec a = (a_1,\dots,a_k)$ in $A$ and $n \ge 1$, consider the \emph{$n$-th matrix numerical range}
\[
  K_{\vec a, n} := \bigl\{(\Phi(a_1),\dots,\Phi(a_k)) : \Phi : A \to M_n(\mathbb{C}) \ \text{ucp}\bigr\},
\]
a compact convex subset of the finite-dimensional space $E_{k,n} = (M_n(\mathbb{C}))^k$; for $n=1$
this is the classical joint numerical range. The key technical point is that $K_{\vec a,n}$ is
computably compact, uniformly in all of its data. This is proved by identifying, via Choi's
correspondence, the support function of $K_{\vec a,n}$ with a constrained supremum over the state
space of $M_n(A)$, and then dualising: the constrained supremum is shown to equal an infimum of a
computable quantity over a computably bounded set of
self-adjoint matrices (Lemma \ref{lem:duality}). With computable compactness in hand, one can decide whether a prescribed tuple $\vec b$ of matrices is approximately attainable, and then
produce an actual computable ucp map attaining it:
 
\begin{theoremB}[Theorem \ref{thm:arveson}]
There is an algorithm which, on input a unital index for a computable presentation $A^\#$ of a unital
\cstar-algebra $A$, integers $n,k \ge 1$, codes for generated points $\vec a$ of $A^\#$, a code for a
rational $\vec b \in E_{k,n}$, and a rational $\delta>0$ such that
$d(\vec b, K_{\vec a,n}) < \delta$, outputs an index for a computable ucp map
$\Phi : A \to M_n(\mathbb{C})$ with $\max_{1 \le i \le k}\|\Phi(a_i) - b_i\| < \delta$.
\end{theoremB}
 
The apparent uncomputability of the Arveson extension theorem is what obstructs the naive attempts
to prove Theorem A, and Theorem B removes that obstruction in exactly the approximate,
matrix-valued form that the definition of nuclearity requires. A pleasant byproduct of the case
$n=1$ is that the computable states of any computable presentation of a unital \cstar-algebra are
weak$^*$-dense in its state space (Corollary \ref{cor:states}); in particular, every such presentation has a
computable state, a fact we use repeatedly and which was previously available to us only through
ad hoc constructions in concrete examples.
 
Sections \ref{sec:categ} and \ref{sec:effkirch} turn to the Cuntz algebra $\mathcal{O}_2$ and to the strongest of the effective
nuclearity conditions we consider. In Section \ref{sec:categ} we prove:
 
\begin{theoremC}[Theorem \ref{thm:categorical}]
For every $n\geq 2$, the Cuntz algebra $\mathcal{O}_n$ is computably categorical, that is, any two computable
presentations of $\mathcal{O}_n$ are computably isomorphic.
\end{theoremC}
 
The proof combines an effective Baire category theorem with the fact that every computable
presentation of $\mathcal{O}_n$ contains a computable \emph{generating} Cuntz $n$-tuple; the
genericity input is a density theorem of R{\o}rdam \cite{Rordam} on unitaries inducing inner automorphisms.
Theorem C allows us to speak of \emph{the} computable presentation of $\mathcal{O}_2$, and it yields
an effective form of the Kirchberg-Phillips absorption theorem: for a simple, separable \cstar-algebra $A$ with
computable presentation $A^\#$, the algebra $A$ is nuclear if and only if the tensor product
presentation $(A \otimes \mathcal{O}_2)^\#$ is computably isomorphic to $\mathcal{O}_2$
(Corollary \ref{cor:effKPnuc}).
 
Kirchberg's theorem also characterises nuclearity of a separable unital \cstar-algebra by the
existence of a unital embedding into $\mathcal{O}_2$ whose image is the range of a conditional
expectation. Effectivizing this characterisation seems more challenging than those conditions considered in Theorem A, and so we
isolate it as Definition \ref{def:effkirchnuc}: $A^\#$ is \emph{effectively Kirchberg nuclear} if there is a
computable embedding $A^\# \hookrightarrow \mathcal{O}_2$ admitting a computable conditional
expectation $\mathcal{O}_2 \to A^\#$. Unlike the conditions in Theorem A, we do not know whether
this one is automatic:
 
\begin{questionA}[Question \ref{question:effkirch}]
Is every computable presentation of a nuclear \cstar-algebra effectively Kirchberg nuclear?
\end{questionA}
 
The difficulty is that the classical proof of Kirchberg's theorem is of a quite different character
from the approximation arguments underlying Theorem A and we see no route to effectivizing it
directly. We offer a useful reformulation of this property, namely that $A^\#$ is effectively
Kirchberg nuclear if and only if it embeds computably, with computable expectation, into
\emph{some} computable presentation of a \emph{simple} nuclear \cstar-algebra (Lemma \ref{easyeffkirchberg}(3)).  We
then exhibit such an ambient algebra in several classes of interest. In the commutative case, the
ambient algebra is the CAR algebra, reached through an effective version of the
Miljutin--Pe{\l}czy\'nski theorem \cite{Miljutin1966,Pelczynski1968} on embeddings of $C(X)$ into $C(2^{\mathbb{N}})$ with
expectation; in the AF case it is a UHF algebra built from a computable AF-certificate.
 
\begin{theoremD}[Corollary \ref{cor:abeliankirch}, Proposition \ref{prop:AFkirch}]
Every computable presentation of a commutative unital \cstar-algebra is effectively Kirchberg
nuclear, and so is every computable presentation of an AF algebra.
\end{theoremD}
 
Since commutative \cstar-algebras other than $\mathbb{C}$ are never simple, Theorem D is not a
consequence of the (straightforward) fact that computable presentations of simple nuclear algebras
are effectively Kirchberg nuclear.
 
For group \cstar-algebras, our approach is dynamical and the point of entry is the observation that
the conditional expectation one needs comes from an invariant state. Precisely, if
$G^\# \curvearrowright^\alpha A^\#$ is a computable action and $\varphi$ is an $\alpha$-invariant state
on $A$ which is computable on $A^\#$, then $\sum_s a_s u_s \mapsto \sum_s \varphi(a_s)u_s$ is a
computable conditional expectation from $A\rtimes_\alpha G$ onto $C^*(G)$ (Lemma \ref{lem:exp}).  (In fact, the existence of \emph{some} conditional expectation onto $C^*(G)$ forces the
existence of an invariant state, so the hypothesis is necessary and sufficient.) Invariant states form a
co-c.e. closed set in the (computably compact) state space of a computably presented nuclear \cstar-algebra (Lemma \ref{invariantstatepi}), and amenability of $G$
guarantees that this class is nonempty; when it is moreover a \emph{singleton}, its unique member is computable
member. This yields:
 
\begin{theoremE}[Theorem \ref{thm:main}, Corollaries \ref{cor:cantor} and \ref{cor:rf}]
Let $G^\#$ be a computable presentation of an amenable group $G$, and suppose $G$ admits a
computable action on a computably compact presentation of a compact metric space which is topologically
free, minimal, and uniquely ergodic. Then $C^*(G)^\#$ is effectively Kirchberg nuclear. In
particular, this holds whenever $G$ admits a computable descending chain of finite-index normal
subgroups with trivial intersection, and hence for every finitely presented, residually finite,
amenable group.
\end{theoremE}
 
Section \ref{sec:completeness} concerns the complexity of nuclearity as a property of indices. Writing $\mathrm{Nuc}$
for the set of indices of computable presentations of nuclear \cstar-algebras, the proof of
Theorem A shows in particular that $\mathrm{Nuc}$ is $\Pi^0_2$ (Proposition \ref{prop:nucispi2}); we conjecture that
it is $\Pi^0_2$-complete (Question \ref{completequestion}). Towards this end, we settle the analogous question for groups,
in the sharper form that the hardness is visible already within the class of (indices of) groups:
 
\begin{theoremF}[Theorems \ref{thm:amenability} and \ref{prop:inneramenable}]
The set $\mathrm{Amen}$ of indices of computable amenable groups is $\Pi^0_2$-complete within the
class $\mathrm{Grp}$ of indices of computable groups. The same is true of the set
$\mathrm{InnAmen}$ of indices of computable inner amenable groups.
\end{theoremF}
 
The construction interleaves two strategies: at stages where a fixed c.e.\ set $W_e$ receives no new
element, we build one more step of a copy of the free group $F_2$ (respectively, of
$SL_2(\overline{\mathbb{F}_2}(t))$), and at stages where $W_e$ grows, we collapse the fragment built
so far into a finite group. That the second strategy is always available rests on a theorem of
Iwasawa \cite{Iwasawa} to the effect that the finite-index normal subgroups of a free group intersect
trivially. 

Transferring this through the universal group \cstar-algebra construction gives the
hardness half of the picture for \cstar-algebras:
 
\begin{theoremG}[Corollary \ref{cor:nuchard}]
The set of indices of c.e.\ presentations of nuclear \cstar-algebras is $\Pi^0_2$-hard.
\end{theoremG}
 
The gap between Theorem G and a complete answer to Question \ref{completequestion} is precisely the question of whether
that set is $\Pi^0_2$ (Question \ref{ques:nuc'pi2}); the proof that $\mathrm{Nuc}$ is $\Pi^0_2$ uses computability
of the presentation in an essential way, through the effective Arveson extension theorem, and so
does not apply to merely c.e.\ presentations.
 
The paper is organised as follows. Section \ref{sec:arveson} develops the effective Arveson extension theorem and
the computable compactness of matrix numerical ranges. Section \ref{sec:effnuc} introduces the effective forms of
nuclearity and proves Theorem A. Section \ref{sec:categ} proves the computable categoricity of $\mathcal{O}_n$
and deduces an effective Kirchberg--Phillips absorption result. Section \ref{sec:effkirch} treats effective Kirchberg
nuclearity, with subsections on commutative algebras, AF algebras, and group \cstar-algebras.
Section \ref{sec:completeness} is devoted to the $\Pi^0_2$ complexity results. 

We assume familiarity with the basics of
computable presentations of metric structures as it pertains to \cstar-algebras; our earlier papers \cite{AFPaper} and \cite{UHFPaper} as well as the recent textbook \cite{EMST} are suitable references.  A good source for the requisite \cstar-algebra theory is \cite{BrownOzawa}.

\subsection*{Conventions}
\emph{Throughout the paper}, all \cstar-algebras are assumed to be unital and all presentations of \cstar-algebras are assumed to be \emph{computably unital}, meaning that the unit is a computable point of the presentation.  Note that a computably unital presentation need not be computable.  A \emph{unital index} for a computable and computably unital presentation consists of an index for the presentation together with an index for the unit.  Since we will be assuming presentations are computably unital, when we say a presentation is computable we mean it is computable and computably unital.

\section*{AI Declaration Statement}

While working on an initial draft of the paper, we uploaded our results to Claude Opus 5 for its opinion.  It suggested an improvement to our results via the Effective Arveson Extension Theorem that is proved in Section 2.  The arguments were reworked and rewritten to accommodate this suggestion.  Also, it suggested an extension and application of our work on computability of crossed products (Proposition \ref{crossedproductcomputable}) which we then incorporated into Section 5.  All of the writing of this paper was done by the authors and as such we take full responsibility for the contents of the paper.

\section{Effective Arveson Extension}\label{sec:arveson}

Throughout this section, we set $E_{k,n} := (\Mn)^k$, which we view as a real vector space equipped with the real inner product
\[
  \langle \vb, \vc\rangle := \Rea \sum_{i=1}^{k}\Tr(c_i^* b_i)
\]
and associated norm $|\cdot|$.  Recall that for a matrix $b \in M_n(\mathbb{C})$, the Hilbert-Schmidt norm $\norm{b}_{\HS} = \Tr(b^*b)^{1/2}$ satisfies $\|b\|_{\mathrm{op}} \le \|b\|_{\HS}$, so bounds measured in $|\cdot|$ give bounds in operator norm coordinatewise.  We view $E_{k,n}$ as equipped with its obvious standard presentation (see also \cite[Subsection 1.2]{UHFPaper}).

We write $\UCP(A,\Mn)$ for the set of ucp maps $A \to \Mn$. Given a tuple
$\va = (a_1,\dots,a_k)$ in $A$, set
\[
  K_{\va,n} := \bigl\{(\Phi(a_1),\dots,\Phi(a_k)) : \Phi \in \UCP(A,\Mn)\bigr\} \subseteq E_{k,n},
\]
which we think of as the \emph{$n$-th matrix (joint) numerical range} of $\va$. For $n = 1$ this is exactly the usual joint numerical
range 
\[W(\vec a):=\{(\omega(a_1),\dots,\omega(a_k)) : \omega \in S(A)\},\] where $S(A)$ denotes the state space of $A$.  The following fact is well-known:

\begin{fact}\label{lem:compact}
$\UCP(A,\Mn)$ is compact in the topology of pointwise weak convergence, and $K_{\va,n}$ is a nonempty
compact convex subset of $E_{k,n}$ satisfying $$|\vb| \le R(\va) := \bigl(n\sum_{i\le k}\|a_i\|^2\bigr)^{1/2}$$
for all $\vb \in K_{\va,n}$.
\end{fact}

Our first main goal is to show that the function $\vec b\mapsto d(\vec b,K_{\va,n}):E_{n,k}\to \mathbb{R}$ is a computable function, uniformly in $n,k,\vec a$.  We will achieve this by showing that the support function for $K_{a,n}$ is a computable function.  (Recall that, for a nonempty, closed, convex set $A\subseteq \mathbb{R}^m$, the support function of $A$ is the function $h_A:K\to \mathbb{R}$ given by $h_A(x):=\sup\{\langle x,a\rangle \ : \ a\in A\}$.)  We first need some preparation.



\subsection{Choi's correspondence with a trace constraint}


Set $\tr_n := \frac1n\Tr$, the normalised trace on $\Mn$. Put
\[
  S_n(A) := \bigl\{\omega \in S(M_n(A)) : \omega\restriction (\Mn \otimes 1) = \tr_n \bigr\}.
\]

\begin{lemma}\label{lem:choi}
The assignments
\[
  \omega_\Phi\bigl((a_{pq})\bigr) := \frac1n \sum_{p,q=1}^{n} \bigl\langle \Phi(a_{pq})\varepsilon_q, \varepsilon_p\bigr\rangle,
  \qquad
  \Phi_\omega(a) := n\sum_{p,q=1}^n \omega(e_{pq}\otimes a)\, e_{pq}
\]
are mutually inverse affine bijections between $\UCP(A,\Mn)$ and $S_n(A)$. (Here $(\varepsilon_p)$ is the
standard basis of $\mathbb{C}^n$.)
\end{lemma}

\begin{proof}
This is the usual Choi correspondence (see, for example, \cite[Proposition 1.5.12]{BrownOzawa}) together with a bookkeeping of the unitality constraint; we record the
computation because we use the exact normalisations.

Let $\Phi$ be cp and let $\xi := \sum_{q}\varepsilon_q \otimes \varepsilon_q \in \mathbb{C}^n\otimes\mathbb{C}^n$,
so $\|\xi\|^2 = n$. Regarding the $n$-th amplification $\Phi_n : M_n(A) \to M_n(\Mn) = B(\mathbb{C}^n \otimes
\mathbb{C}^n)$ (the first tensor factor indexing the block), we have
\[
  \bigl\langle \Phi_n\bigl((a_{pq})\bigr)\xi, \xi\bigr\rangle
  = \sum_{p,q}\bigl\langle \Phi(a_{pq})\varepsilon_q,\varepsilon_p\bigr\rangle
  = n\,\omega_\Phi\bigl((a_{pq})\bigr),
\]
so $\omega_\Phi = \frac1n\langle \Phi_n(\cdot)\xi,\xi\rangle$ is positive; conversely positivity of
$\omega_\Phi$ forces complete positivity of $\Phi$ by the usual argument (apply positivity to
$x^*(a\otimes e_{11})x$-type elements). Taking $(a_{pq}) = e_{pq}\otimes a$
gives $\omega_\Phi(e_{pq}\otimes a) = \frac1n \langle \Phi(a)\varepsilon_q,\varepsilon_p\rangle$, which is
exactly the statement that $\Phi_{\omega_\Phi} = \Phi$; symmetrically $\omega_{\Phi_\omega} = \omega$.

For the constraint: if $\Phi(1) = 1_n$, then $\omega_\Phi(m\otimes 1) = \frac1n\sum_{p,q}m_{pq}\delta_{pq}
= \tr_n(m)$. Conversely, if $\omega\restriction(\Mn\otimes 1) = \tr_n$, then
$\Phi_\omega(1) = n\sum_{p,q}\tr_n(e_{pq})e_{pq} = \sum_{p}e_{pp} = 1_n$. Finally, any positive
$\omega$ with $\omega\restriction(\Mn\otimes1) = \tr_n$ has $\omega(1_B) = \tr_n(1_n) = 1$ and so is a
state.
\end{proof}

For $\vc \in E_{k,n}$ and $\va$ in $A$ put
\[
  y_{\va,\vc} := \tfrac12\bigl(z + z^*\bigr) \in M_n(A)_{\operatorname{sa}},
  \qquad z := \sum_{i=1}^{k} \overline{c_i}\otimes a_i,
\]
where $\overline{c}$ denotes the entrywise complex conjugate of $c$ and
$c \otimes a := (c_{pq}a)_{pq}$.

\begin{lemma}\label{lem:pairing}
For every $\Phi \in \UCP(A,\Mn)$ with associated state $\omega = \omega_\Phi$ and every $\vc \in E_{k,n}$, we have
\[
  \bigl\langle (\Phi(a_1),\dots,\Phi(a_k)), \vc \bigr\rangle = n\,\omega\bigl(y_{\va,\vc}\bigr).
\]
Consequently, letting $h_{\vec a,n}:=h_{K_{\vec a,n}}$ denote the support function of $K_{\vec a,n}$, we have
\begin{equation}\label{eq:support}
  h_{\va,n}(\vc) = n \cdot \sup\bigl\{ \omega(y_{\va,\vc}) : \omega \in S_n(A)\bigr\}.
\end{equation}
Moreover $y_{\va,\vc}$ is a generated point of $M_n(A)^\#$ whenever $\va$ consists of generated points of
$A^\#$ and $\vc$ is rational, and a code for it is computable from codes for $\va, \vc$.
\end{lemma}

\begin{proof}
Writing $\Phi(a)_{qp} = \langle \Phi(a)\varepsilon_p,\varepsilon_q\rangle = n\,\omega(e_{qp}\otimes a)$
by Lemma~\ref{lem:choi}, we have
\[
  \Tr(c_i^*\Phi(a_i)) = \sum_{p,q}\overline{(c_i)_{qp}}\,\Phi(a_i)_{qp}
  = n\,\omega\Bigl(\sum_{p,q}\overline{(c_i)_{qp}}\, e_{qp}\otimes a_i\Bigr)
  = n\,\omega\bigl(\overline{c_i}\otimes a_i\bigr).
\]
Summing over $i$ gives $\sum_i \Tr(c_i^*\Phi(a_i)) = n\,\omega(z)$; since $\omega$ is self-adjoint,
$\Rea\,\omega(z) = \omega(\frac12(z+z^*))$. The last assertion is immediate from the description of the
generated points of $M_n(A)^\#$.
\end{proof}

\subsection{The duality identity}

The point of \eqref{eq:support} is that the constrained supremum on the right has an exactly dual
description as an \emph{infimum over a finite-dimensional space}, and that both the objective and the
constraint region of that infimum are computable:

\begin{lemma}\label{lem:duality}
Let $B$ be a unital \cstar-algebra containing $\Mn$ unitally.  For $z \in B_{\operatorname{sa}}$, set $\lmax(z) := \max \spec(z) =
\inf\{t \in \mathbb{R} : z \le t1\}$.  Set
$S := \{\omega \in S(B) : \omega\restriction \Mn = \tr_n\}$. Then for every $y \in B_{\operatorname{sa}}$, we have
\begin{equation}\label{eq:duality}
  \max_{\omega \in S} \omega(y)
  \;=\;
  \inf_{m \in (\Mn)_{\operatorname{sa}}}\ \bigl[\lmax(y + m) - \tr_n(m)\bigr].
\end{equation}
Moreover, the infimum on the right hand side of the above equation is in fact attained.
\end{lemma}

First, we will need the following:

\begin{lemma}\label{lem:coercive}
In the situation of Lemma~\ref{lem:duality}, put $g(m) := \lmax(y+m)-\tr_n(m)$ for $m \in (\Mn)_{\operatorname{sa}}$.
Then:
\begin{enumerate}
\item $g(m + t1_n) = g(m)$ for all $t \in \mathbb{R}$.  In particular, the infimum in (\ref{eq:duality}) may be taken over the
trace $0$ self-adjoint matrices.
\item $|g(m) - g(m')| \le 2\|m-m'\|$ for all $m,m'\in (\Mn)_{\operatorname{sa}}$.
\item If $m\in (\Mn)_{\operatorname{sa}}$ has trace $0$, then $\lmax(m) \ge \|m\|/n$, and thus
$$g(m) \ge \|m\|/n - \|y\|.$$
\item Consequently, the right hand side of (\ref{eq:duality}) equals
\[
  \min\bigl\{ \lmax(y+m) \;:\; m \in (\Mn)_{\operatorname{sa}},\ \Tr(m)=0,\ \|m\| \le 2n\|y\| \bigr\}.
\]
\end{enumerate}
\end{lemma}

\begin{proof}
(1) follows from the equations $\lmax(y+m+t1) = \lmax(y+m)+t$ and $\tr_n(m+t1_n) = \tr_n(m)+t$.

(2) From $z \le z' + \|z-z'\|1$ we get $\lmax(z)\le\lmax(z') + \|z-z'\|$ and hence
$|\lmax(z)-\lmax(z')|\le\|z-z'\|$; also $|\tr_n(m)-\tr_n(m')|\le\|m-m'\|$.

(3) Let $\mu_1 \ge \cdots \ge \mu_n$ be the eigenvalues of $m$, so $\sum_p \mu_p = 0$ and
$\|m\| = \max(\mu_1, -\mu_n)$. If $\|m\| = \mu_1$ we are done. Otherwise $\|m\| = -\mu_n$ and $n \ge 2$
(a traceless self-adjoint $1\times1$ matrix is $0$), and $\sum_{p<n}\mu_p = -\mu_n = \|m\|$, so
$\mu_1 \ge \|m\|/(n-1) \ge \|m\|/n$. For the second inequality, use the inequality
$\lmax(y+m) \ge \lmax(m) - \|y\|$, which follows from the more general inequality mentioned in the proof of item (2).

(4) By (1), we may take the infimum over $m$ of trace $0$. Next note that $g(0) = \lmax(y) \le \|y\|$.  If, on the other hand, $m$ has trace $0$ and 
$\|m\| > 2n\|y\|$, we have $g(m) > 2\|y\| - \|y\| = \|y\| \ge g(0)$ by (3). So the infimum over $m$ of trace $0$ equals the infimum over the indicated (compact) ball, which is therefore attained.
\end{proof}


\begin{proof}[Proof of Lemma \ref{lem:duality}]
Write $\alpha$ for the left-hand side; note that the supremum is attained since $S$ is weak$^*$-compact and
$S \ne \emptyset$ since $\tr_n$ extends to a state of $B$ by Hahn--Banach, positivity being automatic
for a norm-one extension of a state on a unital subalgebra.  Let $\beta$ denote the right-hand side
infimum, which we know equals the minimum appearing in item (4) of the previous lemma.

We first show that $\alpha \le \beta$). For $\omega \in S$ and $m \in (\Mn)_{\operatorname{sa}}$, we have
\[
  \omega(y) = \omega(y+m) - \omega(m) = \omega(y+m) - \tr_n(m) \le \lmax(y+m) - \tr_n(m),
\]
using that $\omega$ is a state and $y + m$ is self-adjoint.

We now show that $\beta \le \alpha$. Suppose, towards a contradiction, that this is not the case and fix $\gamma \in \mathbb{R}$ with
$\alpha < \gamma < \beta$. Set $V := (\Mn)_{\operatorname{sa}} \subseteq B_{\operatorname{sa}}$ (again identified with its image).

We first dispose of a degenerate case: if $y \in V$, say $y = m_0$, then every $\omega \in S$ has
$\omega(y) = \tr_n(m_0)$, so $\alpha = \tr_n(m_0)$; taking $m = -m_0$ in the right-hand side of
\eqref{eq:duality} gives $\beta \le \lmax(0) + \tr_n(m_0) = \alpha$, contradicting $\alpha<\beta$. So
$y \notin V$ and $W := \mathbb{R}y \oplus V$ is a direct sum. Define a real-linear functional $\varphi$ on $W$ by
\[
  \varphi(sy + m) := s\gamma + \tr_n(m) \qquad (s \in \mathbb{R},\ m \in (\Mn)_{\operatorname{sa}}).
\]

\noindent Claim: $\varphi \le \lmax$ on $W$. 

\

\noindent To prove the claim, we consider three cases:
\begin{itemize}
\item $s = 0$: $\tr_n$ is a state of $\Mn$ and $\lmax(m)$ is computed in $\Mn$ as well (a unital
inclusion of \cstar-algebras preserves spectra), so $\tr_n(m) \le \lmax(m)$.
\item $s > 0$: dividing by $s$ and writing $m' = m/s$, the claim reads
$\gamma \le \lmax(y+m') - \tr_n(m')$, which holds for every $m'$ because $\gamma < \beta$.
\item $s < 0$: dividing by $|s|$ and writing $m' = m/|s|$, the claim reads
$-\gamma + \tr_n(m') \le \lmax(-y+m')$. Pick $\omega \in S$ with $\omega(y) = \alpha$; then
$\lmax(-y+m') \ge \omega(-y+m') = -\alpha + \tr_n(m') > -\gamma + \tr_n(m')$ since $\alpha<\gamma$.
\end{itemize}

Now $\lmax : B_{\operatorname{sa}} \to \mathbb{R}$ is sublinear, so by the
Hahn-Banach theorem there is a real-linear $\tilde\varphi : B_{\operatorname{sa}}\to
\mathbb{R}$ extending $\varphi$ with $\tilde\varphi \le \lmax$. Then $\tilde\varphi$ is positive: if
$z \ge 0$ then $\lmax(-z) \le 0$, so $-\tilde\varphi(z) = \tilde\varphi(-z)\le 0$. Also
$\tilde\varphi(1) = \varphi(1_n) = \tr_n(1_n) = 1$. Complexifying,
$\omega_0(z) := \tilde\varphi(\Rea z) + i\tilde\varphi(\operatorname{Im} z)$ is a state of $B$; it satisfies
$\omega_0 \restriction \Mn = \tr_n$ (true on self-adjoints by construction, hence everywhere by
linearity), so $\omega_0 \in S$, and $\omega_0(y) = \varphi(y) = \gamma > \alpha$, contradicting the
definition of $\alpha$.

Finally, the infimum $\beta$ is attained by
Lemma~\ref{lem:coercive}.
\end{proof}

\subsection{Computability of the support function}

Using the lemma from the previous subsection, we can now prove:

\begin{proposition}\label{prop:support}
Let $A^\#$ be a computable presentation of a \cstar-algebra $A$, let $\vec a=(a_1,\ldots,a_k)$ be a $k$-tuple of generated points of $A$, and let $n\geq 1$.  Then the support function $h_{\vec a,n}$ of $K_{\vec a,n}$ is computable.  Moreover, an index for $h_{\vec a,n}$ can be computed from a unital index for $A^\#$, indices for $a_1,\ldots,a_k$, and $n$.
\end{proposition}

\begin{proof}
By Lemmas~\ref{lem:pairing}, \ref{lem:duality} and Lemma~\ref{lem:coercive}(4) applied to $B = M_n(A)$, $y := y_{\va,\vc}$, we have $h_{\va,n}(\vc) = n\beta$, where
\[
  \beta = \min\bigl\{\lmax(y+m\otimes1) : m \in (\Mn)_{sa},\ \Tr(m)=0,\ \|m\|\le 2n\|y\|\bigr\}.
\]
We check that each ingredient is computable.

First, since $M_n(A)^\#$ is a computable presentation uniformly in $n$ (see \cite[Prop.~1.7]{UHFPaper}), a code for the generated point $y$ of $M_n(A)^\#$ is computable from
$\va,\vc$, and, moreover, $\|y\|$ is computable . Fix a rational $N > \|y\| + 2n\|y\| + 1$.

Next, for $m$ a rational self-adjoint matrix, $y + m\otimes 1$ is a computable point of
$M_n(A)^\#$ with a modulus computable from $m$, and
$\|y+m\otimes1\| < N$. Since $y + m\otimes 1 + N1 \ge 0$, we have
\[
  \lmax(y+m\otimes1) = \|y + m\otimes 1 + N\cdot 1\| - N,
\]
and the right-hand side is computable uniformly in $m$.

Finally, the set
$C := \{m \in (\Mn)_{\operatorname{sa}} : \Tr(m) = 0,\ \|m\|\le 2n\|y\|\}$ is a compact subset of a
$(n^2-1)$-dimensional real vector space with a computable rational basis, and from a rational upper
bound for $\|y\|$ we can compute a finite rational $\delta$-net of $C$ (of size computable in $\delta$)
lying in the slightly larger set $$C' := \{m : \Tr(m)=0,\ \|m\|\le 2n\|y\|+1\};$$ note that, by
Lemma~\ref{lem:coercive}(4), the minimum of $\lmax(y+\cdot)$ over $C'$ still equals $\beta$. By
Lemma~\ref{lem:coercive}(2) the function $\lmax(y+\cdot)$ is $2$-Lipschitz in $m$, so evaluating it to precision
$\varepsilon/(4n)$ on an $(\varepsilon/(8n))$-net of $C'$ and taking the minimum returns a rational
within $\varepsilon/n$ of $\beta$, hence within $\varepsilon$ of $n\beta$.
\end{proof}

\begin{corollary}\label{cor:hausdorff}
Uniformly in the data of Proposition~\ref{prop:support}, we have:
\begin{enumerate}
\item The distance function $\vb \mapsto d(\vb, K_{\va,n}):E_{k,n}\to \mathbb{R}$ is computable.
\item $K_{\va,n}$ is a computably compact subset of $E_{k,n}$.
\end{enumerate}
\end{corollary}

\begin{proof}
(1) It is a well-known fact that, for a nonempty compact convex $K \subseteq E$ with support function $h$, one has
\begin{equation}\label{eq:distformula}
  d(x, K) = \max_{|c|\le 1}\bigl(\langle x, c\rangle - h(c)\bigr).
\end{equation}
By Proposition~\ref{prop:support}, the function being maximized is computable.  Consequently, its
maximum over the (computably compact) unit ball is computable.

(2) follows from general computable analysis:  in $M_n(\C)$, closed rational balls are computably compact; in such a space, a compact set satisfying (1) is automatically computably compact \cite[Proposition 3.29]{Downey.Melnikov.2023}.
\end{proof}

\subsection{The effective Arveson extension theorem}

We can now prove the main result of this section:

\begin{theorem}[Effective $\varepsilon$-Arveson]\label{thm:arveson}
There is an algorithm that, given
\begin{itemize}
\item a unital index for a computable presentation $A^\#$ of a unital \cstar-algebra $A$,
\item $n,k \ge 1$ and codes for generated points $\va = (a_1,\dots,a_k)$ of $A^\#$, and
\item a code for a rational $\vb \in E_{k,n}$ and a rational $\delta > 0$ such that
$d(\vb, K_{\va,n}) < \delta$,
\end{itemize}
outputs an index for a computable ucp map $\Phi : A \to \Mn$ with
\[
  \max_{1\le i \le k}\|\Phi(a_i) - b_i\| < \delta .
\]
\end{theorem}

\begin{proof}
Fix an effective enumeration $x_1, x_2, x_3, \dots$ of the generated points of $A^\#$ with $x_i = a_i$
for $i \le k$. For $l \ge k$ set:
\begin{itemize}
    \item $\vx^{(l)} := (x_1,\dots,x_l)$;
    \item $K_l := K_{\vx^{(l)},n} \subseteq
E_{l,n}$, and
\item $\pi_l : E_{l+1,n}\to E_{l,n}$ for the projection onto the first $l$ coordinates.
\end{itemize}

Since $K_l$ and $K_{l+1}$ are the images of the \emph{same} set $\UCP(A,\Mn)$ under evaluation at
$\vx^{(l)}$ and $\vx^{(l+1)}$ respectively, we have
\begin{equation}\label{eq:onto}
  \pi_l(K_{l+1}) = K_l .
\end{equation}

Using Corollary~\ref{cor:hausdorff}, compute a rational $\delta_0$ with
$d(\vb,K_k) < \delta_0 < \delta$, and fix a rational $\theta>0$ with $\delta_0 + 4\theta < \delta$. Set
\[
  \eta_k := \delta_0, \qquad \eta_l := \theta\, 2^{-(l-k)} \quad (l > k).
\]

We construct rational $v^{(l)} \in E_{l,n}$ for $l \ge k$ with
\begin{equation}\label{eq:induct}
  d\bigl(v^{(l)}, K_l\bigr) < \eta_l
  \qquad\text{and}\qquad
  \bigl|\pi_l\bigl(v^{(l+1)}\bigr) - v^{(l)}\bigr| < \eta_l + \eta_{l+1}.
\end{equation}
Start with $v^{(k)} := \vb$, which satisfies the first condition. Given $v^{(l)}$, by Corollary~\ref{cor:hausdorff}(1)) we may search for a rational
$v^{(l+1)}$ satisfying the above two conditions; the search terminates because such a point exists. Indeed, choose
$u \in K_l$ with $|u - v^{(l)}| < \eta_l$, lift it by \eqref{eq:onto} to $u' \in K_{l+1}$ with
$\pi_l(u') = u$, and take any rational $v^{(l+1)}$ with $|v^{(l+1)} - u'| < \eta_{l+1}$; then
$d(v^{(l+1)},K_{l+1}) < \eta_{l+1}$ and
$|\pi_l(v^{(l+1)}) - v^{(l)}| \le |\pi_l(v^{(l+1)}) - u| + |u - v^{(l)}| < \eta_{l+1}+\eta_l$.

By \eqref{eq:induct}, for each $i$ the sequence
$\bigl(v^{(l)}_i\bigr)_{l \ge \max(i,k)}$ satisfies
$|v^{(l+1)}_i - v^{(l)}_i| \le \eta_l + \eta_{l+1}$, and $\sum_{l\ge k}(\eta_l+\eta_{l+1})$ converges
with a computable rate. As a result,
\[
  \Phi(x_i) := \lim_{l\to\infty} v^{(l)}_i
\]
exists and is a computable matrix, uniformly in $i$. Moreover
for $i \le k$, we have
\[
  |\Phi(a_i) - b_i| \le \sum_{l \ge k}(\eta_l + \eta_{l+1})
  \le \eta_k + 2\sum_{l>k}\eta_l = \delta_0 + 4\theta < \delta,
\]
and $\|\cdot\|_{\mathrm{op}} \le \|\cdot\|_{\HS}$ gives the stated norm estimate.

We next note that $\Phi$ is the restriction of a ucp map $A\to \Mn$. To see this, fix $l$. For $l' > l$, \eqref{eq:induct} gives
$|\pi_l(v^{(l')}) - v^{(l)}| \le \sum_{j=l}^{l'-1}(\eta_j+\eta_{j+1})$, so letting $l'\to\infty$, we have
\[
  d\Bigl(\bigl(\Phi(x_1),\dots,\Phi(x_{l'})\bigr), K_{l'}\Bigr) \le \eta_{l'} + \sum_{j\ge l'}(\eta_j+\eta_{j+1})
  \;\xrightarrow[l'\to\infty]{}\; 0 .
\]
Since $\pi_l(K_{l'}) = K_l$ by \eqref{eq:onto} and projections are contractive, the left-hand side
dominates $d\bigl((\Phi(x_1),\dots,\Phi(x_l)), K_l\bigr)$, which is therefore $0$; as $K_l$ is closed, we have
\[
  \bigl(\Phi(x_1),\dots,\Phi(x_l)\bigr) \in K_l \qquad \text{for every } l .
\]
Thus, for each $l$ there is $\Phi_l \in \UCP(A,\Mn)$ with $\Phi_l(x_i) = \Phi(x_i)$ for $i \le l$. By
Fact~\ref{lem:compact}, the net $(\Phi_l)$ has a point-weak cluster point
$\Phi_\infty \in \UCP(A,\Mn)$, and for each fixed $i$ we have $\Phi_l(x_i) = \Phi(x_i)$ for all
$l \ge i$, so $\Phi_\infty(x_i) = \Phi(x_i)$. Thus $\Phi_\infty$ is a ucp map agreeing with our
computed values on all generated points.  Being contractive, it is the unique continuous extension of
those values, and the algorithm above witnesses that it is computable.
\end{proof}

The previous theorem has the following interesting consequence:

\begin{corollary}\label{cor:states}
Let $A^\#$ be a computable presentation of a unital \cstar-algebra. Then, uniformly in a unital index for $A^\#$, we have:
\begin{enumerate}
\item The joint numerical range of any tuple of generated points is a computably compact convex subset
of $\mathbb{C}^k$.
\item Given generated points $\va$, a rational $\vb \in \mathbb{C}^k$, and rational $\delta>0$ with
$d(\vb, W(\va)) < \delta$, one can compute a state $\omega$ of $A$, computable as a functional on
$A^\#$, with $\max_i|\omega(a_i) - b_i| < \delta$.
\item In particular, the states on $A$ that are computable with respect to $A^\#$ are weak$^*$-dense in
$S(A)$.
\end{enumerate}
\end{corollary}

\begin{proof}
(1) and (2) follow by taking $n = 1$ in Corollary~\ref{cor:hausdorff} and Theorem~\ref{thm:arveson}, using that
$\UCP(A,\mathbb{C}) = S(A)$ and $S_1(A) = S(A)$.  (Note that, in this case, \eqref{eq:duality} degenerates to the triviality
$\max_{\omega\in S(A)}\omega(y) = \lmax(y)$, and \eqref{eq:support} says the support function of $W(\va)$ is
computable because norms are.)  (3) is an immediate consequence of (2).
\end{proof}

\section{Effective nuclearity}\label{sec:effnuc}

\begin{definition}\label{def effnuc}
Suppose that $A^\#$ is a presentation of a \cstar-algebra $A$.  We say that $A^\#$ is:
\begin{enumerate}
\item \textbf{effectively nuclear} if there is an algorithm such that, upon input a finite tuple $\vec a$ of generated points of $A^\#$ and rational $\epsilon>0$, returns:
    \begin{itemize}
        \item a positive integer $n\in \mathbb{N}$, and
        \item indices for computable ucp maps $\phi:A^\#\to M_n(\mathbb{C})$ and $\psi:M_n(\mathbb{C})\to A^\#$
    \end{itemize}
such that $$\|(\psi\circ \phi)(\vec a)-\vec a\|<\epsilon.$$
\item \textbf{effectively Choi-Effros nuclear} if there is a uniformly computable sequence $\phi_n:A^\#\to A^\#$ of finite rank ucp maps such that:
\begin{itemize}
    \item there is a computable function $f:\mathbb{N}\to \mathbb{N}$ such that $\operatorname{rank}(\phi_n)\leq f(n)$ for all $n\in \mathbb{N}$, and
    \item $\phi_n$ computably converges pointwise to the identity, that is, there is an algorithm such that, given a generated point $a$ of $A^\#$ and $\epsilon>0$, returns $n\in \mathbb{N}$ such that $\|\phi_m(a)-a\|<\epsilon$ for all $m\geq n$.
\end{itemize} 
\end{enumerate}
\end{definition}

Recall our standing assumption is that the presentations of \cstar-algebras that we consider must have the unit as a computable point.  In fact, this property follows from being effectively nuclear, as does computability of the presentation.

\begin{proposition}\label{lem:automaticallycomputable}
Suppose that $A^\#$ is a presentation (not assumed to be computable or computably unital) of a unital \cstar-algebra $A$.  If $A^\#$ is effectively nuclear then $A^\#$ is computable and computably unital, and from an index for the algorithm witnessing the effective nuclearity of $A^\#$ we can compute a unital index for $A^\#$.
\end{proposition}
\begin{proof}
Fix a generated point $a$ of $A^\#$ and $\epsilon>0$.  Compute ucp maps $\phi:A\to \Mn$ and $\psi:\Mn\to A$ such that $\|(\psi\circ \phi)(a)-a\|<\epsilon$.  Note then that $\|a\|\geq \|\phi(a)\|\geq \|\psi(\phi(a))\|\geq \|a\|-\epsilon$.  Since $\phi$ is computable, we can effectively approximate $\|\phi(a)\|$ to within $\epsilon$, resulting in an effective approximation of $\|a\|$ within $2\epsilon$.  Moreover, the standard presentation of $\Mn$ allows us to compute an index for the $n \times n$ identity matrix $1_n$, and thus we can compute an $A^\#$-index for computable point $1_A$ since $1_A = \psi(1_n)$.
\end{proof}

\begin{question}
If $A^\#$ is effective Choi-Effros nuclear, does it follows that $A^\#$ is computable and computably unital?
\end{question}

We return to our convention that presentations of \cstar-algebras are computably unital.

Before proving the main theorem of this section, we note the following result, which simply notes that the fact that perturbing an almost unital cp map to a ucp map can be carried out effectively:

\begin{lemma}\label{perturb}
There is an an algorithm such that, upon input a unital index for a  presentation $A^\#$ of a \cstar-algebra $A$ and a computable cpc map $\psi':M_n(\mathbb{C})\to A$ satisfying $\|\psi'(1)-1\|<\epsilon$ with $\epsilon\leq 1/2$, returns the index for a computable ucp map $\psi:\Mn\to A$ such that $\|\psi-\psi'\|\leq 2\epsilon$.
\end{lemma}

\begin{proof}
Since $\epsilon<1$, we have that $\psi'(1)$ is an invertible, positive element of $A$.  We thus know that $\psi'(1)^{-1/2}$ exists.  The desired ucp map is defined by $$\psi(x):=\psi'(1)^{-1/2}\psi'(x)\psi'(1)^{-1/2}.$$  
One checks easily that $\|\psi-\psi'\|\leq \epsilon/(1-\epsilon)\leq 2\epsilon$.  If $\psi'$ is computable, then $\psi'(1)$ is computable, whence so is $\psi'(1)^{-1/2}$ by the computable functional calculus \cite[Corollary 4.6]{EagleMcNicholl}.  
\end{proof}

We again recall that our presentations are assumed to be computably unital.
\begin{theorem}\label{mainequivalence}
Suppose that $A^\#$ is a presentation of a \cstar-algebra $A$.  The following are equivalent:
\begin{enumerate}
    \item $A^\#$ is computable and $A$ is nuclear.
    \item $A^\#$ is effectively nuclear.
    \item $A^\#$ is computable and effectively Choi-Effros nuclear.
\end{enumerate}
\end{theorem}

\begin{proof}
The implications (2)$\Rightarrow$ (3) and (3)$\Rightarrow$ (1) are clear.  To see that (1) implies (2), suppose one is given a finite $k$-tuple $\vec a$ of generated points of $A^\#$ and a rational $\epsilon>0$.  Choose rational $\delta\in (0,\min(\epsilon/4,1/2))$.  Search for a triple $(n,\vec b,w)$, where $n\geq 1$, $\vec b$ is a rational $k$-tuple from $M_n(\mathbb{C})$, and $w$ is a generated point of $M_n(A)^\#$ such that:
\begin{itemize}
\item $d(\vec b,K_{\vec a,n})<\delta$,
    \item $\|\sum_{jj}\pi_{jj}(w^*w)-1\|<\delta$
    \item $\max_{1\leq i\leq k}\|a_i-\sum_{jl}\pi_{jl}(b_i)\pi_{jl}(w^*w)\|<\delta$.
\end{itemize}

By Corollary \ref{cor:hausdorff} above, this can be effectively done, and since $A$ is nuclear, the search is successful.  By Theorem \ref{thm:arveson} above, we can then compute an index for a computable ucp map $\phi:A\to M_n(\mathbb{C})$ such that $\max_{1\leq i\leq k}\|\phi(a_i)-b_i\|<\delta$.  Let $\psi':M_n(\mathbb{C})\to A$ be the unique linear map for which $\psi'(e_{jl})=\pi_{jl}(w^*w)$.  Then $\psi'$ is a computable cpc map with $\|\psi'(1)-1\|<\delta$.  Note that 
$$\|\psi'(\phi(a_i))-a_i\|\leq \|\psi'(\phi(a_i))-\psi'(b_i)\|+\|\psi'(b_i)-a_i\|\leq 2\delta.$$
By Lemma \ref{perturb}, we may computably perturb $\psi'$ to a computable ucp map $\psi:M_n(\mathbb{C})\to A$ such that $\|\psi-\psi'\|<2\delta$.  It follows that $\|\psi(\phi(a_i))-a_i\|<4\delta<\epsilon$.   
\end{proof}
We note that the proof of Theorem \ref{mainequivalence} is uniform if we are given a unital index for $A^\#$.  On the other hand, it cannot be made uniform given just an index for $A^\#$ without the index for the unit.  Indeed, the third author has shown in \cite[Theorem 5.3]{mcnicholl2024evaluative} that there is a commutative (hence nuclear) unital \cstar-algebra $B$ for which there is no algorithm taking each index for a computable presentation of $B$ to an index for the unit in that presentation, although the unit is a computable point in every computable presentation of $B$ by \cite[Theorem 4.1]{burton2024computable}.  By Proposition \ref{lem:automaticallycomputable}, for this $B$ there is no algorithm taking each index for a computable presentation of $B$ to an index for the algorithm witnessing effective nuclearity.

\section{Computable categoricity of $\mathcal{O}_2$}\label{sec:categ}

For separable, \emph{simple} \cstar-algebras $A$, the Kirchberg-Phillips theorem states that $A$ is nuclear if and only if $A\otimes \mathcal{O}_2\cong \mathcal{O}_2$.  It is tempting to effectivize this condition by insisting that $(A\otimes \mathcal{O}_2)^\#$ be computably isomorphic to $\mathcal{O}_2$, where $\mathcal{O}_2$ is equipped with its standard presentation (which is computable, see \cite[Corollary 3.16]{fox2022computable}) and $(A\otimes \mathcal{O}_2)^\#$ is the associated tensor product presentation (see \cite[Subsection 2.3]{Goldbring.2024+}).  In this section, we will see that such a condition is equivalent to $A^\#$ being computable and nuclear.  
This will follow from a result that is important in its own right, namely that $\mathcal{O}_2$ is \emph{computably categorical}, meaning any two computable presentations of $\mathcal{O}_2$ are computably isomorphic.  We now work towards proving that goal.

\begin{definition}
A tuple $(r_1, \ldots, r_n)$ in $\mathcal{O}_n$ is a \emph{Cuntz $n$-tuple} if it satisfies the Cuntz relations
\begin{itemize}
    \item $r_j^*r_j = 1$ for all $j=1,\ldots,n$
    \item $\sum_{j=1}^nr_jr_j^* = 1$
\end{itemize}
A \emph{generating Cuntz $n$-tuple} is a Cuntz $n$-tuple that generates the entire ambient algebra $\mathcal{O}_n$.
\end{definition}

The proof of computable categoricity of $\mathcal{O}_n$ relies on the fact that every computable presentation of $\mathcal{O}_n$ contains a computable generating Cuntz $n$-tuple, that is, a generating Cuntz $n$-tuple consisting of computable points of the presentation.  To prove this, we will need the following effective version of the Baire category theorem:

\begin{lemma}\label{lemma:DenseGDelta}
Let $X^\#$ be a computable presentation of a metric space.  Suppose that $D_0 \supseteq D_1 \supseteq \cdots$ are c.e. open sets of $X^\#$, each of which is dense in $X$.  Then the set of computable points in $\bigcap_{n \in \mathbb{N}}D_n$ is dense in $X$.
\end{lemma}
\begin{proof}
Fix a rational open ball $B$ of $X^\#$.  Since $D_0$ is dense in $X$, $B \cap D_0 \neq \emptyset$.  Since $B$ is a rational ball and $D_0$ is c.e. open, we can find a rational open ball $B_0$ such that $B_0 \subseteq B \cap D_0$ and $\operatorname{diam}(B_0) < 1$.  Now $B_0 \cap D_1 \neq \emptyset$, and $B_0$ is a rational open ball and $D_1$ is c.e. open, so we can find a rational open ball $B_1$ such that $B_1 \subseteq B_0 \cap D_1$ and $\operatorname{diam}(B_1) < 1/2$.  Continuing in this way we obtain rational open balls $B \supseteq B_0 \supseteq B_1 \supseteq \cdots$ with each $B_{j+1} \subseteq B_j \cap D_{j+1}$ and $\operatorname{diam}(B_j) < 2^{-j}$.  The centres of the balls $B_0, B_1, \ldots$ thus form a computably Cauchy sequence that converges to a computable point in $B \cap \bigcap_{n \in \mathbb{N}}D_n$.
\end{proof}

\begin{lemma}\label{lemma:ComputableGeneratingCuntz}
For every $n \in \mathbb{N}$ and every computable presentation $\mathcal{O}_n^\#$ of $\mathcal{O}_n$, there is a computable generating Cuntz $n$-tuple of $\mathcal{O}_n^\#$.
\end{lemma}
\begin{proof}
The Cuntz relations are computably weakly stable (see \cite[Example 12]{fox2024locally}), so by \cite[Theorem 1.12]{UHFPaper} the set of Cuntz $n$-tuples is a c.e. closed set of $\mathcal{O}_n^\#$.  In particular, we can find a computable Cuntz $n$-tuple, say $(r_1, \ldots, r_n)$.  For any unitary $u$, $(ur_1, \ldots, ur_n)$ is clearly a Cuntz $n$-tuple.  It therefore suffices to find a computable unitary $u$ such that $(ur_1, \ldots, ur_n)$ is generating.

Let $U(\mathcal{O}_n)$ denote the set of unitaries in $\mathcal{O}_n$.  Let $(a_k)_{k \in \mathbb{N}}$ be a computable enumeration of the generated points of $\mathcal{O}_n^\#$, and let $(p_k(x_1, \ldots, x_n))_{k \in \mathbb{N}}$ be a computable enumeration of the rational *-polynomials in $n$ variables.  For each $m, k \in \mathbb{N}$ and each rational $\epsilon > 0$, define
\[D_{m,k,\epsilon} = \{u \in U(\mathcal{O}_n) : \norm{p_m(ur_1, \ldots, ur_n) - a_k} < \epsilon\}.\]
For each $k, \epsilon$, let $D_{k, \epsilon} = \bigcup_{m \in \mathbb{N}}D_{m, k, \epsilon}$, and let $D = \bigcap_{k \in \mathbb{N}, \epsilon \in \mathbb{Q}^{<0}}D_{k,\epsilon}$.  Observe that $u \in D$ if and only if $(ur_1, \ldots, ur_n)$ is a generating Cuntz $n$-tuple.  It thus remains to be shown that there is a computable point in $D$.

Let $(s_1, \ldots, s_n)$ be a generating Cuntz $n$-tuple.  Let
\[Y = \{u \in U(\mathcal{O}_n) : (us_1, \ldots, us_n) \text{ is generating}\}.\]
Note that $u \in Y$ if and only if the endomorphism determined by $(s_1, \ldots, s_n) \mapsto (us_1, \ldots, us_n)$ is an automorphism.  R\o rdam \cite{Rordam} showed that the set of unitaries for which the corresponding endomorphism is an \emph{inner} automorphism is dense (see \cite[Proposition 7]{BratteliKishimoto} for a self-contained proof).  It follows that $Y$ is also dense.  Let $v$ be the (unique) unitary such that $(r_1, \ldots, r_n) = (vs_1, \ldots, vs_n)$.  Then for any unitary $u$ we have $u \in D$ if and only if $uv \in Y$, so $D = Yv^*$.  Since multiplication on the right by $v^*$ is an isometry of $U(\mathcal{O}_n)$, it follows that $D$ is dense in $U(\mathcal{O}_n)$, and hence also that each $D_{k, \epsilon}$ is dense in $U(\mathcal{O}_n)$.

The property of being a unitary is given by computably weakly stable relations, so $U(\mathcal{O}_n)$ is a c.e. closed set of $\mathcal{O}_n^\#$.  In particular, it contains a dense set of computable points, which we take to be the computable presentation $U(\mathcal{O}_n)^\#$ of $U(\mathcal{O}_n)$ induced by $\mathcal{O}_n^\#$.  Observe that each $D_{m, k, \epsilon}$ is a c.e. open set of $U(\mathcal{O}_n)^\#$, and therefore each $D_{k, \epsilon}$ is also a c.e. open set of $U(\mathcal{O}_n)^\#$.  Thus Lemma \ref{lemma:DenseGDelta} implies that there is a computable point in $D$.
\end{proof}

With the previous lemma in hand, we can now prove the main result of this section:

\begin{theorem}\label{thm:categorical}
For every $n\geq 2$, the Cuntz algebra $\mathcal{O}_n$ is computably categorical.
\end{theorem}
\begin{proof}
Let $\mathcal{O}_n^\#$ be a computable presentation of $\mathcal{O}_n$.  Let $\mathcal{O}_n^\dagger$ be the standard presentation of $\mathcal{O}_n$, that is, the presentation whose special points are a generating Cuntz $n$-tuple $(s_1, \ldots, s_n)$.  By \cite[Corollary 3.16]{fox2022computable}, $\mathcal{O}_n^\dagger$ is a computable presentation.  By Lemma \ref{lemma:ComputableGeneratingCuntz} there is a generating Cuntz $n$-tuple $(r_1, \ldots, r_n)$ with each $r_j$ a computable point of $\mathcal{O}_n^\#$.  The map induced by mapping each $r_j$ to $s_j$ is thus a computable isomorphism from $\mathcal{O}_n^\#$ onto $\mathcal{O}_n^\dagger$.
\end{proof}

\begin{corollary}\label{cor:effKPnuc}
If $A$ is a separable, simple \cstar-algebra and $A^\#$ is a computable presentation of $A$, then $A$ is nuclear if and only if there is a computable isomorphism $(A\otimes \mathcal{O}_2)^\#\cong \mathcal{O}_2$.
\end{corollary}

We remark that the proof of Theorem \ref{thm:categorical} does not seem to cover the case of the Cuntz algebra $\mathcal{O}_\infty$ (whose standard presentation is computable by \cite[Corollary 3.16]{fox2022computable}).

\begin{question}
    Is the Cuntz algebra $\mathcal{O}_\infty$ computably categorical?
\end{question}
\section{Effectively Kirchberg nuclear}\label{sec:effkirch}

\subsection{General discussion}

\begin{definition}\label{def:effkirchnuc}
   We say that $A^\#$ is \textbf{effectively Kirchberg nuclear} if there is a computable embedding $A^\#\hookrightarrow \mathcal O_2$ which admits a computable conditional expectation map $\mathcal O_2\to A^\#$.
\end{definition}

Note that $A^\#$ being effectively Kirchberg nuclear implies that $A^\#$ is computable; in fact, this follows only from the assumption that there is a computable embedding $A^\#\hookrightarrow \mathcal{O}_2$.

\begin{lemma}\label{easyeffkirchberg}

\

\begin{enumerate}
\item Suppose that there is a computable embedding $A^\#\hookrightarrow B^\#$ that admits a computable expectation $B^\#\to A^\#$.  Further suppose that $B^\#$ is effectively Kirchberg nuclear.  Then $A^\#$ is effectively Kirchberg nuclear.
\item If $A^\#$ is a computable presentation of a \emph{simple} nuclear \cstar-algebra $A$, then $A^\#$ is effectively Kirchberg nuclear.
    \item For a presentation $A^\#$ of a \cstar-algebra $A$, we have that $A^\#$ is effectively Kirchberg nuclear if and only if there is a computable presentation $B^\#$ of a simple nuclear \cstar-algebra $B$ and a computable embedding $A^\#\hookrightarrow B^\#$ with computable expectation $B^\#\to A^\#$.
    \end{enumerate}
\end{lemma}

\begin{proof}
(1) is immediate from the definitions.  For (2), note that the obvious embedding $A^\#\hookrightarrow (A\otimes \mathcal{O}_2)^\#$ is computable.  Moreover, if $\phi$ is a computable state on $\mathcal{O}_2$ (see Corollary \ref{cor:states}), the slice map $a\otimes b\mapsto \phi(b)a:A\otimes \mathcal{O}_2\to A$ is a computable conditional expectation for this embedding.  It remains to use (1) and computable categoricity of $\mathcal{O}_2$.  (3) then follows immediately from (1) and (2).
\end{proof}

\begin{remark}
In the above proof, one can explicitly construct a computable state on $\mathcal{O}_2$.  Indeed, if we view the generators of $\mathcal O_2$ as the transformations $T_1$ and $T_2$ on the separable Hilbert space $H$ with orthonormal basis $(e_n)_{n\in \mathbb{N}}$ for which $T_1(e_n)=e_{2n}$ and $T_2(e_n)=e_{2n+1}$, then the vector state $\langle \cdot e_0,e_0\rangle$ is easily seen to be computable.
\end{remark}

\begin{question}\label{question:effkirch}
Is every computable presentation of a nuclear \cstar-algebra effectively Kirchberg nuclear?
\end{question}

While we do not have an answer to this question at the moment, in the next sections we give an affirmative answer for all computable presentations of commutative \cstar-algebras, AF \cstar-algebras, and group \cstar-algebras associated to certain (but not all) amenable groups.

\subsection{Commutative algebras}
In this subsection we show that computable presentations of commutative unital \cstar-algebras are effectively Kirchberg nuclear.  Aside from $\mathbb{C}$, commutative \cstar-algebras are never simple, so Lemma \ref{easyeffkirchberg}(2) does not apply.

We begin by establishing effective Kirchberg nuclearity for the \cstar-algebra $C(2^{\mathbb{N}})$ of continuous complex-valued functions on the Cantor space.  This \cstar-algebra is computably categorical \cite[Corollary 6.7]{mcnicholl2024evaluative}, so we do not need to specify which computable presentation to consider.

\begin{theorem}\label{thm:CantorKirchberg}
$C(2^\mathbb{N})$ is effectively Kirchberg nuclear.
\end{theorem}
\begin{proof}
Let $\boldB$ be the CAR algebra, that is, the UHF algebra with supernatural number $2^{\infty}$.  The CAR algebra is the universal \cstar-algebra generated by a sequence $(b_n)_{n \in \mathbb{N}}$ satisfying the CAR relations for all $n$ and $m$:
\[b_nb_m+b_mb_n = 0 \qquad \text{ and } \qquad b_nb_m^*+b_m^*b_n = \begin{cases}1 & n=m \\ 0 & n\neq m\end{cases}.\]
Let $\boldB^\#$ be the presentation of $\boldB$ whose special points are $(b_n)_{n \in \mathbb{N}}$.  By \cite[Theorem 3.11]{fox2022computable}, the presentation $\boldB^\#$ is c.e., and since $\boldB$ is simple, by \cite[Theorem 3.14]{fox2022computable}, $\boldB^\#$ is computable.  For each $n$, let $d_n = b_n^*b_n$.  Let $\boldD$ be the \cstar-subalgebra of $\boldB$ generated by $(d_n)_{n \in \mathbb{N}}$.  It is well-known that $\boldD \cong C(2^{\mathbb{N}})$.  Let $\boldD^\#$ be the presentation of $\boldD$ whose special points are also $(d_n)_{n \in \mathbb{N}}$.  There is a uniform procedure for writing any generated point of $\boldD^\#$ as a generated point of $\boldB^\#$, whence computability of $\boldD^\#$ follows from computability of $\boldB^\#$.  As mentioned above, $C(2^{\mathbb{N}})$ is computably categorical; by Lemma \ref{easyeffkirchberg}(3), it thus suffices to show that the computable embedding of $\boldD^\#$ into $\boldB^\#$ admits a computable conditional expectation.

Indeed, there is a canonical conditional expectation map for this embedding, defined as follows.  Fix $n \in \mathbb{N}$ and $t = (t_1, \ldots, t_n) \in \{0, 1\}^n$, and define
\[P_t := \prod_{i=1}^nd_i^{t_i}(1-d_i)^{1-t_i}.\]
The expectation map $E : \boldB \to \boldD$ is defined on elements of $C^*(b_1, \ldots, b_n)$ by
\[E(x) := \sum_{t \in \{0, 1\}^n}P_txP_t.\]
In particular, this shows that $E$ is computable on generated points of $\boldB^\#$.  Since $E$ is contractive, this suffices to show that $E$ is a computable map from $\boldB^\#$ to $\boldD^\#$.
\end{proof}

Miljutin \cite{Miljutin1966} showed (in different terminology) that $C([0, 1])$ embeds into $C(2^\mathbb{N})$ with a conditional expectation.  This result was then extended from $[0, 1]$ to arbitrary compact metric spaces by Pe\l czy\'nski \cite{Pelczynski1968}.  We now prove an effective version of this result.  

Suppose that $X$ is a compact Polish space.  Recall from \cite{mcnicholl2024evaluative} that a presentation $C(X)^\dagger$ of $C(X)$ is \emph{evaluative over} a presentation $X^\#$ of $X$ if the map $(f, p) \mapsto f(p)$ is from the presentation of $C(X) \times X$ induced by $C(X)^\dagger$ and $X^\#$.  It was shown in \cite[Theorem 1.3]{mcnicholl2024evaluative} that, given a computably compact presentation $X^\#$ of $X$, there is a unique (up to computable isomorphism) computable presentation $C(X^\#)$ of $C(X)$ that is evaluative over $X^\#$.

\begin{lemma}\label{lemma:EComputable}
Suppose that $X^\#$ is a computably compact presentation of a compact Polish space $X$.  For each $f : X \to \mathbb{C}$, the following are equivalent:
\begin{itemize}
    \item $f$ is a computable map with respect to $X^\#$.
    \item $f$ is a computable point of $C(X^\#)$.
\end{itemize}
Moreover, a $C(X^\#)$-name for $f$ can be computed from an index for $f$ as a computable map, and vice-versa.
\end{lemma}
\begin{proof}
Suppose first that $f$ is a computable map.  Given any computable $g : X^\# \to \mathbb{C}$, the function $|f-g|$ is also computable, and hence by \cite[Proposition 3.5(1)]{Downey.Melnikov.2023}, $\norm{f-g}$ is computable, uniformly in the pair $(f, g)$.  Let $(a_n)_{n \in \mathbb{N}}$ be an effective enumeration of the generated points of $C(X^\#)$.  By definition of $C(X^\#)$, each $a_n$ is a computable map from $X^\#$ to $\mathbb{C}$, uniformly in $n$.  Therefore, given $k \in \mathbb{N}$, we can search until an $n(k)$ is located such that $\norm{f-a_{n(k)}} < 2^{-k}$.  Then $(a_{n(k)})_{k\in\mathbb{N}}$ is a $C(X^\#)$-name for $f$, so $f$ is a computable point of $C(X^\#)$.

Now suppose that $f$ is a computable point of $C(X^\#)$.  Let $E : C(X) \times X \to \mathbb{C}$ be the evaluation map $E(f,x) = f(x)$.  Since $C(X^\#)$ is evaluative over $X^\#$, $E$ is computable.  It follows that $f$ is a computable map, and an index for $f$ can be computed from a $C(X^\#)$-name for $f$ and an index for $E$.
\end{proof}

\begin{theorem}\label{thm:effMilyutin}
Let $X^\#$ be a computably compact presentation of a compact Polish space $X$.  Then there is a computable embedding of $C(X^\#)$ into $C(2^{\mathbb{N}})$ with a computable conditional expectation.
\end{theorem}
\begin{proof}
Using the computable compactness of $X^\#$, we can compute a positive integer $N$ such that $\operatorname{diam}(X) < N$.  We next use computable compactness to construct a tree of open balls in $X$.  Specifically, we compute: 
\begin{itemize}
    \item An increasing sequence $(k_s)_{s \in \N}$ of positive integers. 

    \item An array $(B_{s,n})_{s \in \N, n < k_s}$ of 
    rational open balls of $X^\#$. 

    \item An array $(j_{s + 1, n})_{s \in \N, n \leq k_s}$ of nonnegative integers.
\end{itemize}
We require that:
\begin{enumerate}
    \item The formal diameter of $B_{s,n}$ is at most $2^{-s} N$. 

    \item $X = \bigcup_{n < k_s} B_{s,n}$. 

    \item $0 = j_{s + 1, 0} < \ldots < j_{s+1, k_s} = k_{s + 1}$. 

    \item $\overline{B_{s+1, n}} \subseteq B_{s,n'}$ when 
    $j_{s+1, n'} \leq n < j_{s + 1, n' + 1}$. 
\end{enumerate}

Call a sequence $\sigma \in \N^{< \omega} \setminus \{\emptyset\}$ \emph{valid} if
$\sigma(0) < k_0$ and $j_{s + 1, \sigma(s)} \leq \sigma(s + 1) < j_{s + 1, \sigma(s) + 1}$ whenever $s + 1 < |\sigma|$.  Observe that a non-empty subsequence of a valid sequence is again valid.
Let 
\[T = \{\sigma  \in \mathbb{N}^{<\omega}\setminus\{\emptyset\} : \sigma \text{ is valid}\} \cup \{\emptyset\}.\]

When $\sigma \in T$, set 
\[
B_\sigma := \left\{
\begin{array}{cc}
X & \sigma = \emptyset\\
B_{|\sigma| - 1, \sigma(|\sigma| - 1)} & \sigma \neq \emptyset \\
\end{array}
\right. 
\]
Let $[T]$ denote the set of branches through $T$, equipped with the usual  ultrametric $d(\beta_1, \beta_2) = 2^{-n}$, where $n$ is the length of the longest common initial segment of $\beta_1$ and $\beta_2$.  Note that $T$ is a computable tree with no dead ends, and it is computably branching, in the sense that there is a computable binary function $B$ such that if $\sigma \in T$ and $\rho \in T$ extends $\sigma$, then for all $i < |\rho|$, $\rho(i) \leq B(\sigma, i)$.  It follows that there is a computable presentation $[T]^\#$ of $[T]$ and a computable surjection from $2^{\mathbb{N}}$ onto $[T]$ (see \cite[Theorem 3.40, proof of Claim 1]{Downey.Melnikov.2023}; the presentation has special points of the form $v_\sigma$ where $\sigma \in T$ and $v_\sigma$ is the leftmost branch of $T$ extending $\sigma$); by \cite[Lemma 3.31]{Downey.Melnikov.2023}, $[T]^\#$ is computably compact.

Since $[T]$ is a compact $0$-dimensional Hausdorff space without isolated points, it is homeomorphic to $2^{\mathbb{N}}$, and hence $C([T]) \cong C(2^{\mathbb{N}})$.  By computable categoricity of $C(2^{\mathbb{N}})$, $C([T]^\#)$ is computably isomorphic to the standard presentation of $C(2^{\mathbb{N}})$.  Thus it suffices to find a computable embedding of $C(X^\#)$ into $C([T]^\#)$ with a computable conditional expectation.

Conditions (1) and (4) ensure that if $\beta \in [T]$, then $\{\overline{B_\sigma} : \sigma \subset \beta\}$ is a nested family of closed balls whose diameters shrink to $0$, and moreover $\bigcap_{\sigma \subset \beta}B_\sigma = \bigcap_{\sigma \subset \beta}\overline{B_\sigma}$. We define $F : [T] \to X$ by setting $F(\beta)$ to be the unique point of $\bigcap_{\sigma \subset \beta}B_\sigma$.  It follows from condition (2) that $F$ is surjective.  Condition (1) ensures that the centers of the balls $B_\sigma$ for $\sigma \subset \beta$ form a computably Cauchy sequence in $X$, and hence $F$ is a computable map from $[T]^\#$ to $X^\#$.   Let $\tilde{F} : C(X) \to C([T])$ be the Gelfand dual of $F$.  By \cite[Proposition 6.2]{mcnicholl2024evaluative}, $\tilde{F}$ is a computable embedding of $C(X^\#)$ into $C([T]^\#)$.  We next turn to finding the associated computable conditional expectation map.

By \cite[Theorem 3]{Yasgui.1999}, for each $s$ we can find a computable partition of unity subordinate to the cover $\{B_{s,n} : n < k_s\}$, and moreover we can do this uniformly in $s$.  More specifically, there exists an array $(\psi_{s,n})_{s \in \N, n < k_s}$
so that:
\begin{enumerate}
    \item $\psi_{s,n}$ is a computable map of $X^\#$ into $[0,1]$ uniformly in 
    $s,n$.
    \item $\psi_{s,n}(x) = 0$ whenever $x \not\in B_{s,n}$.
    \item For each $s$, $\sum_{n} \psi_{s,n} = 1$.
\end{enumerate}

For each $p \in X$, define $\mu_p$ to be the unique Borel probability measure on 
$[T]$ so that, for each $\sigma \in T$, we have
\[
\mu_p([\sigma]) = \prod_{j < |\sigma|} \psi_{j, \sigma(j)}(p)
\]
(where we define the product over the empty set to be $1$). 

Let $E : C([T]) \to C(X)$ be defined by
\[E(g)(p) = \int_{F^{-1}(\{p\})}g\,d\mu_p.\]
This map $E$ is a conditional expectation for the embedding $\tilde{F} : C(X) \to C([T])$.  
  Our construction ensures that for each $p$ we have $\mu_p([T] \setminus F^{-1}(\{p\})) = 0$, so we can equivalently write
\[E(g)(p) = \int_{[T]} g\,d\mu_p.\]
It remains to be shown that $E$ is computable as a map from $C([T]^\#)$ to $C(X^\#)$.  Let $S : C([T]) \times X \to \mathbb{C}$ be defined by
\[S(g,p) = \int_{[T]}g\,d\mu_p.\]

Let $P([T])$ denote the space of Radon probability measures on $[T]$ with the Prokhorov metric (which we denote by $d_P$). 
Define $P([T]^\#)$ to be the presentation of $P([T])$ whose special points are the convex rational linear combinations of point masses at the special points of $[T]^\#$.  It was shown in \cite[Proposition 4.1.1]{HoyrupRojas} that the presentation $P([T]^\#)$ is computable.

Suppose we have computable points $g$ of $C([T]^\#)$ and $p$ of $X^\#$.  We describe an algorithm for computing $S(g, p)$.  First, since the map $(q, \sigma) \mapsto \mu_q([\sigma])$ is computable, it follows that $q \mapsto \mu_q$ is a computable map from $X^\#$ to $P([T]^\#)$.  Thus we can compute a sequence $(\mu_n)_{n \in \mathbb{N}}$ in $P([T])$ of special points of $P([T]^\#)$ that effectively converges to $\mu_p$.  By definition of the presentation $P([T]^\#)$ and the fact that $C([T]^\#)$ is evaluative over $T^\#$, the map $n \mapsto \int_{[T]}g\,d\mu_n$ is computable, uniformly in $g$.  

By \cite[Theorem 4.3 and Theorem 4.10]{Rojas}, $(\mu_n)_{n \in \mathbb{N}}$ uniformly effectively weakly converges to $\mu_p$.  By definition, this means that from a name for $f \in C([T])$ and a bound on $\norm{f}$ we can compute a modulus of convergence for $(\int_{[T]}f\,d\mu_n)_{n\in\mathbb{N}}$ converging to $\int_{[T]}f\,d\mu_p$.  Since $g$ is a computable point of $C([T]^\#)$ we can compute $\norm{g}$.  We therefore have a computable sequence of computable complex numbers that effectively converges to $S(g,p)$, uniformly in $g$ and $p$, as required to show that $S$ is computable.  Since $S$ is computable, the map $p \mapsto S(g,p)$ is computable, uniformly in $g$.  That is, $E(g) : [T]^\# \to \mathbb{C}$ is a computable map, uniformly in $g$.  By Lemma \ref{lemma:EComputable} it follows that $E(g)$ is a computable point of $C([T]^\#)$, uniformly in $g$, so $E$ is computable.
\end{proof}

Lemma \ref{easyeffkirchberg} and Theorems \ref{thm:CantorKirchberg} and \ref{thm:effMilyutin} immediately imply:

\begin{corollary}\label{cor:abeliankirch}
Every computable presentation of a commutative unital \cstar-algebra is effectively Kirchberg nuclear.
\end{corollary}
\subsection{AF-algebras}

In this subsection, we establish that any computable presentation of an AF-algebra is effectively Kirchberg nuclear.

\begin{proposition}\label{prop:AFkirch}
Suppose that $A$ is an AF algebra with computable presentation $A^\#$.  Then $A^\#$ is effectively Kirchberg nuclear.
\end{proposition}

\begin{proof}
By Lemma \ref{easyeffkirchberg}, it suffices to show that $A^\#$ admits a computable embedding into a computable presentation $B^\#$ of a UHF algebra $B$ that admits a computable expectation.  To carry this out, fix a computable AF-certificate $(F_j,\psi_j)$ for $A^\#$ (see \cite[Main Theorem 1]{AFPaper}).  One may then compute unital embeddings $F_j\hookrightarrow M_{k(j)}(\mathbb{C})$ into full matrix algebras in a coherent way.  Set $B:=\lim M_{k(j)}(\mathbb{C})$, which is then a UHF algebra.  By construction, the embedding $A^\#\hookrightarrow B^\#$ is computable, where $B^\#$ is the inductive limit presentation of $B$ corresponding to the standard presentation the $M_{k(j)}$'s.  Finally, note that the embeddings $F_j\hookrightarrow M_{k(j)}(\mathbb{C})$ admit canonical conditional expectation maps, which are clearly computable, and moreover these conditional expectation maps are coherent throughout the sequence.  It follows that they yield a computable conditional expectation map $B^\#\to A^\#$, completing the proof.
\end{proof}
\subsection{Group \cstar-algebras}

In this section, we show that certain group \cstar-algebras associated to amenable groups are effectively Kirchberg nuclear.  Our approach is dynamical and so we need to take a short detour through crossed product \cstar-algebras and their presentations.

First, we recall that a \emph{presented group} (in the sense of Fox \cite[Definition 3.19]{fox2022computable}) is a pair $G^\#:=(G,(g_n)_{n\in \N})$, where $G$ is a group and $(g_n)_{n\in \N}$ is a countable sequence in $G$ generating $G$ as a group; in this case, we say that $G^\#$ is a \emph{presentation of $G$}.  We say that the presentation $G^\#$ is \emph{computable} if there is an algorithm such that, upon input a finite word in the sequence $(g_n)_{n\in \N}$, determines whether or not the word equals the identity in $G$.

Suppose now that $G\curvearrowright^\alpha A$ is an action of a group $G$ on a \cstar-algebra $A$.  Then the universal crossed product \cstar-algebra $A\rtimes_\alpha G$ can be described by the universal \cstar-algebra generated by unitaries $u_s$ for $s\in G$ and the elements of $A$ subject to the relations on $G$ and the relations stating that $\alpha_s(a)=u_sau_s^*$ for all $s\in G$ and $a\in A$.  (If $G=\langle S|R\rangle$, then the above description can consider only unitaries coming from $S$, relations coming from $R$, and conjugation relations using only elements from $S$.) As such, if $G^\#$ and $A^\#$ are presentations of $G$ and $A$ respectively, one obtains a canonical induced presentation $(A\rtimes_\alpha G)^\#$ of $A\rtimes_\alpha G$.  In this context, we say that the action $G^\#\curvearrowright^\alpha A^\#$ is \textbf{computable} if there is an algorithm such that, upon input of a special point $g$ of $G^\#$, generated point $a$ of $A^\#$, and rational $\epsilon>0$, returns a generated point $b$ of $A^\#$ such that $\|\alpha_g(a)-b\|<\epsilon$.  If this is the case and $G^\#$ and $A^\#$ are additionally assumed to be computable, then $(A\rtimes_\alpha G)^\#$ is c.e.  (In fact, it suffices for $G^\#$ and $A^\#$ to be merely c.e.)

Let $C_c(G,A)$ denote the set of finite formal sums $\sum_{s\in G}a_su_s$ with each $a_s\in A$.  $C_c(G,A)$ is naturally a $*$-algebra under the twisted multiplication and involution operations given by
$$\left(\sum_s a_su_s\right)\left(\sum_sb_su_s\right)=\sum_{s,t}a_s\alpha_s(b_t)u_{st}$$ and
$$\left(\sum_sa_su_s\right)^*=\alpha_{s^{-1}}(a_s^*)u_{s^{-1}}.$$
Note that $A$ can be viewed as a $*$-subalgebra of $C_c(G,A)$ by identifying $a\in A$ with the formal sum $\sum_s a_su_s$ with $a_e=a$ and $a_s=0$ for $s\not=e$.  Note also that $C_c(G,A)$ sits naturally as a dense subalgebra of $A\rtimes_\alpha G$.  

Suppose now that $A\subseteq B(H)$ is a faithful representation of $A$.  We define a representation of the $*$-algebra $C_c(G,A)$ on the Hilbert space $H\otimes \ell^2(G)$ by defining $au_s(v\otimes \delta_t):=\alpha_{t^{-1}}(a)(v)\otimes \delta_{st}$.  One defines the reduced crossed product $A\rtimes_{r,\alpha}G$ to be the \cstar-subalgebra of $B(H\otimes \ell^2(G))$ generated by the image of $C_c(G,A)$.  It is known that this \cstar-algebra is independent on the choice of faithful representation $A\subseteq B(H)$ (see, for example, \cite[Proposition 4.1.5]{BrownOzawa}).  If $G^\#$ and $A^\#$ are presentations of $G$ and $A$ respectively, then there is a canonical induced presentation $(A\rtimes_{r,\alpha} G)^\#$ of $A\rtimes_{r,\alpha} G$ consisting of finite formal sums $\sum_s a_su_s$ with each $a_s$ a generated point of $A^\#$.  

Once again suppose that $A^\#$ is a presentation of $A$ and further assume that $H^\#$ is a presentation of $H$.  We say that the representation $A^\#\subseteq B(H^\#)$ is \textbf{computable} if there is an algorithm such that, upon input a generated point $a$ of $A^\#$, a generated point $v$ of $H^\#$, and $\epsilon>0$, returns a generated point $w$ of $H^*$ such that $\|av-w\|<\epsilon$.  The following lemma follows immediately from the definitions:

\begin{lemma}
Suppose that $G^\#$ is a computable presentation, $G^\#\curvearrowright^\alpha A^\#$ is a computable action, and $A^\#\subseteq B(H^\#)$ is a computable representation.  Then the induced presentation $(A\rtimes_{r,\alpha} G)^\#$ is left-c.e.
\end{lemma}

If the action $G\curvearrowright^\alpha A$ is \textbf{amenable} (for example, when $G$ is amenable), then the canonical surjective morphism $A\rtimes_\alpha G\to A\rtimes_{r,\alpha}G$ is an isomorphism (see \cite[Theorem 4.3.4]{BrownOzawa}).  In this case, we have:

\begin{proposition}\label{crossedproductcomputable}
Suppose that $G^\#$ is a computable presentation of $G$, $A^\#$ is a computable presentation of $A$, and $G^\#\curvearrowright^\alpha A^\#$ is a computable, amenable action.  Furthermore, suppose that there is a computable representation $A^\#\subseteq B(H^\#)$.  Then the induced presentation $(A\rtimes_{\alpha} G)^\#$ is computable.
\end{proposition}

We now need to discuss situations where an embedding of a group \cstar-algebra into a crossed product admits a conditional expectation.  This is essentially standard fare; we are simply paying attention to the computability aspects of the expectation: 

\begin{lemma}\label{lem:exp}
Let $G^\# \curvearrowright^\alpha A^\#$ be a computable action of a computable presentation of a group
on a computable presentation of a unital \cstar-algebra, and let $\varphi$ be a $\alpha$-invariant state
on $A$ which is computable with respect to $A^\#$. Then:
\begin{enumerate}
\item The canonical map $\iota : C^*(G) \to A \rtimes_\alpha G$, $u_s \mapsto 1_A u_s$, is an isometric
$*$-isomorphism onto its image, and it is computable as a map $C^*(G)^\# \to (A\rtimes_\alpha G)^\#$.
\item The formula
\[
  E\Bigl(\sum_{s} a_s u_s\Bigr) := \sum_s \varphi(a_s)\, u_s
\]
defines a conditional expectation $E : A \rtimes_\alpha G \to C^*(G)$ with $E \circ \iota = \mathrm{id}$,
and $E$ is computable as a map $(A\rtimes_\alpha G)^\# \to C^*(G)^\#$.
\end{enumerate}
\end{lemma}

\begin{proof}
That $\iota$ is an isometric $*$-isomorphism onto its image and that $E$ is a conditional expectation map is probably standard fare, but we sketch a proof for the reader's convenience.  Let $(\pi,H,\xi)$ be the result of applying the GNS construction to $\varphi$.  Since $\varphi$ is invariant, we obtain a unitary representation $V$ of $G$ on $H$ so that $V_g\pi(a)V_g^*=\pi(\alpha_g(a))$ for all $a\in A$ and such that $V_g\xi=\xi$.  We then get a $*$-homomorphism $\Lambda:A\rtimes_\alpha G\to B(H)\otimes_{\max}C^*(G)$ given by $\Lambda(au_g)=\pi(a)V_g\otimes u_g$.  The expectation $E$ is the result of composing with the ucp slice map $\omega_\xi\otimes \operatorname{id}$.  Since $E\circ \iota=\operatorname{id}$, we see that $\iota$ is isometric.   

The computability of $\iota$ is immediate.  It remains to see that $E$ is computable.  To see this, given a generated point $x = \sum_{s \in F} a_s u_s$ of
$(A\rtimes_\alpha G)^\#$ with $F \subseteq G$ finite, and rational $\epsilon>0$, compute rationals $q_s$
with $|\varphi(a_s)-q_s| < \epsilon/|F|$; then $\sum_s q_s u_s$ is a generated point of $C^*(G)^\#$ with
\[
  \Bigl\| E(x) - \sum_s q_s u_s \Bigr\| \le \sum_{s\in F} |\varphi(a_s)-q_s| < \epsilon. \qedhere
\]
\end{proof}

In order to apply the previous lemma, we need to know when computable invariant states exist.  In order to achieve that goal, we first mention the following recent fact of the first-named author:

\begin{fact}[\cite{eagle}]
Suppose that $A^\#$ is a computable presentation of a nuclear \cstar-algebra $A$.  Then the state space $S(A)$, with the weak$^*$ topology, admits a computably compact presentation
$S(A)^\#$ such that the evaluation map $(\phi, a) \mapsto \phi(a)$ is computable on the induced presentation of $S(A) \times A$.
\end{fact}

\begin{lemma}\label{invariantstatepi}
Let $A^\#$ be a computable presentation of a nuclear \cstar-algebra and $G^\# \curvearrowright^\alpha A^\#$
a computable action. Then the set $S(A)^G$ of $\alpha$-invariant states is a co-c.e. closed set of $S(A)^\#$.
\end{lemma}

\begin{proof}
Fix $g \in G$, $i \in \omega$, and rational $\epsilon>0$. Using computability of the action,
compute (an index for) a generated point $a_{j(g,i,\epsilon)}$ with
$\|\alpha_g(a_i) - a_{j(g,i,\epsilon)}\| < \epsilon$. We claim that $S(A)^G$ is the set of
$\varphi \in S(A)$ satisfying
\begin{equation}\label{eq:inv}
  \bigl| \varphi(a_{j(g,i,\epsilon)}) - \varphi(a_i) \bigr| \le 2\epsilon
  \qquad\text{for all } g \in G,\ i \in \omega,\ \epsilon \in \mathbb{Q}_{>0},
\end{equation}
which is a countable computable conjunction of $\Pi^0_1$ conditions on computable coordinates, hence
$\Pi^0_1$. Indeed, if $\varphi$ is invariant then
$|\varphi(a_{j}) - \varphi(a_i)| = |\varphi(a_j) - \varphi(\alpha_g(a_i))| \le \epsilon$; conversely if
\eqref{eq:inv} holds then $|\varphi(\alpha_g(a_i)) - \varphi(a_i)| \le 3\epsilon$ for every $\epsilon$,
so $\varphi \circ \alpha_g$ and $\varphi$ agree on generated points and hence, being bounded, on all of
$A$.
\end{proof}

We can now present our main criterion for effective Kirchberg nuclearity of group \cstar-algebras for amenable groups:

\begin{theorem}\label{thm:main}
Let $G^\#$ be a computable presentation of an amenable group, let $A^\#$ be a computable presentation
of a nuclear \cstar-algebra admitting a computable representation, and let $G^\# \curvearrowright^\alpha A^\#$ be a computable action
such that
\begin{enumerate}
\item $A \rtimes_\alpha G$ is simple, and
\item $\alpha$ admits a \emph{unique} invariant state on $A$.
\end{enumerate}
Then $C^*(G)^\#$ is effectively Kirchberg nuclear. Moreover the conclusion is uniform in indices for
$G^\#$, $A^\#$ and $\alpha$.
\end{theorem}

\begin{proof}
Since $G$ is amenable, $C^*(G) = C^*_r(G)$ and, by \cite[Corollary 3.22]{fox2022computable}, $C^*(G)^\#$ is a computable
presentation. By Proposition \ref{crossedproductcomputable},  $(A \rtimes_\alpha G)^\#$ is a computable presentation of $A\rtimes G$.  Moreover, $A\rtimes_\alpha G$ is nuclear,
being the crossed product of a nuclear algebra by an amenable group (see \cite[Theorem 4.3.4]{BrownOzawa}).

By Lemma \ref{invariantstatepi}, $S(A)^G$ is a co-c.e. closed set in the computably compact space $S(A)^\#$, and
it is nonempty by amenability; by hypothesis (2), it is a singleton $\{\varphi\}$.  Therefore, $\varphi$ is a computable point of $S(A)^\#$ and is thus a computable state on $A^\#$.   Lemma~\ref{lem:exp} now supplies a computable
embedding $\iota : C^*(G)^\# \hookrightarrow (A\rtimes_\alpha G)^\#$ together with a computable
conditional expectation $E : (A\rtimes_\alpha G)^\# \to C^*(G)^\#$.

Since $A\rtimes_\alpha G$ is simple and nuclear, Lemma \ref{easyeffkirchberg}(2)
gives that $(A\rtimes_\alpha G)^\#$ is effectively Kirchberg nuclear.  By Lemma \ref{easyeffkirchberg}(1), it follows that $C^*(G)^\#$ is effectively Kirchberg nuclear.
\end{proof}

We now specialize to the case of groups acting on commutative \cstar-algebras.  In order to apply the above theorem, we need to know computable presentations of commutative algebras always admit computable representations:

\begin{lemma}
Suppose $X^\#$ is a computably compact presentation.  Then the associated evaluative presentation $C(X^\#)$ of $C(X)$ admits a computable faithful state, and thus a computable faithful representation.
\end{lemma}

\begin{proof}
Suppose that $(x_i)_{i\in \mathbb{N}}$ are the special points of $X$ and set $\mu=\sum_{i=1}^\infty 2^{-i}\delta_{x_i}$, a probability measure on $X$.  Define $\phi:C(X)\to \mathbb{C}$ by $\phi(f):=\int fd\mu$.  Then $\phi$ is a faithful state (faithful as the measure is strictly positive).  It is also computable since 
$$\left\|\int fd\mu-\sum_{i=1}^{N} 2^{-i}f(x_i)\right\|\leq 2^{-N}\|f\|_\infty$$ and the presentation $C(X^\#)$ is computable and evaluative.  The fact that $C(X^\#)$ admits a computable faithful representation follows from applying the GNS construction with respect to the state $\phi$ (see \cite{eagle} for more details).
\end{proof}

Recall that if $G\curvearrowright^\alpha X$ is an action, then one obtains an induced action $G\curvearrowright^\alpha C(X)$ defined by $\alpha(g,f)(x):=f(\alpha(g^{-1},x))$.

\begin{lemma}\label{computableactionlemma}
Suppose that $G^\#$ is a computable presentation of $G$, $X^\#$ is a computably compact presentation of $X$, and the action $G^\#\curvearrowright^\alpha X^\#$ is a computable action.  Then the induced action $G^\#\curvearrowright^\alpha C(X^\#)$ is also computable.
\end{lemma}

\begin{proof}
By Lemma \ref{lemma:EComputable}, it suffices to show that $\alpha(g,f)$ is a computable function, uniformly in indices for $g$ and $f$.  However, this is immediate from computability of the action and from $C(X^\#)$ being computable and evaluative.
\end{proof}

\begin{corollary}\label{cor:cantor}
Let $G^\#$ be a computable presentation of an amenable group admitting a computable action on a
computably compact presentation $X^\#$ of a metric space $X$ which is
\begin{enumerate}
\item topologically free,
\item minimal, and
\item uniquely ergodic.
\end{enumerate}
Then $C^*(G)^\#$ is effectively Kirchberg nuclear.
\end{corollary}

\begin{proof}
$C(X^\#)$ is a computable presentation of $C(X)$ by \cite[Theorem 4.3]{fox2022computable}, and the induced action on it is computable by Lemma \ref{computableactionlemma}. Topological freeness together with
minimality gives simplicity of $C(X)\rtimes G$ by a result of Archbold and Spielberg \cite{archbold1994topologically}, and unique ergodicity is
precisely the uniqueness of the invariant state. It remains to apply Theorem~\ref{thm:main}.
\end{proof}

\begin{corollary}\label{cor:rf}
Let $G^\#$ be a computable presentation of an amenable group.  Suppose further that there is a computable descending
sequence $N_0 \supseteq N_1 \supseteq \cdots$ of finite-index normal subgroups of $G$ with
$\bigcap_k N_k = \{1\}$, where computable meaning that membership in $N_k$ is decidable uniformly in
$k$ as is the index $[G:N_k]$. Then $C^*(G)^\#$ is effectively
Kirchberg nuclear.
\end{corollary}

\begin{proof}
Let $X := \varprojlim_k G/N_k$.  We first equip $X$ with a computable presentation.  For $g \in G$, write $\hat g := (gN_k)_{k\in\N} \in X$. The map $g \mapsto \hat g:G\to X$ is an injective group
homomorphism with dense image.  
Let $(w_n)_{n\in\N}$ be an effective enumeration of the words in the generators of $G^\#$ and their
inverses. Define $X^\#$ to be the presentation of $(X,d)$ whose $n$-th special point is
$\widehat{w_n}$.  To see that this presentation is computable, observe first that for words $u,v$ and $k \in \N$, we have
\[
  \widehat u_k = \widehat v_k \iff uN_k = vN_k \iff u^{-1}v \in N_k ,
\]
which is decidable uniformly in $(u,v,k)$. Since the $N_k$ are nested, the set of $k$ at which
$\widehat u$ and $\widehat v$ agree is an initial segment of $\N$. Hence the algorithm to compute the metric proceeds by first deciding whether $u^{-1}v = 1$ in $G$. If so, output $d(\widehat u, \widehat v) = 0$.  Otherwise, search for the least $k$ with $u^{-1}v \notin N_k$ and output $2^{-k}$.  The second step indeed terminates whenever $u^{-1}v \ne 1$ as $\bigcap_k N_k = \{1\}$.

We next observe that $X^\#$ is computably compact.  To see this, fix $j$ and compute a transversal $t_{1,j},\ldots,t_{m_j,j}$ of $N_j$.  Then the finite set
$\{\widehat{t_{i,j}} : 1\leq i \leq m_j\}$ of special points is easily seen to be a $2^{-(j+1)}$-net of $X$.

We next note that $G$ acts on $X$ by left translation; the action is computable since membership in each $N_k$ is
decidable. The action is free since $\bigcap_k N_k = \{1\}$ by normality and minimal since $G$ acts transitively on each
finite level $G/N_k$. As the $N_k$ are normal, $X$ is a profinite group and $G$ acts by translation
through a dense image, so the Haar measure of $X$ is the unique invariant measure and unique ergodicity
holds. 

The result now follows from
Corollary~\ref{cor:cantor}.
\end{proof}

\begin{examples}
We list some examples of computable presentations of amenable groups for which the computable chain appearing in the previous theorem exists and thus the reduced group \cstar-algebra is effectively Kirchberg nuclear. 
\begin{enumerate}
\item Suppose that $G$ is a finitely presented, residually finite amenable group.  One may then find such a computable chain of quotients as in the previous theorem.  Indeed, the finite quotients of a
finitely presented group are computably enumerable (only finitely many relators need checking), so one
may search for finite quotients separating a given element from the identity and intersect the kernels
found so far. 
    \item The wreath product $\mathbb{Z}_2 \wr \mathbb{Z}$ is finitely generated but not finitely presented, yet the
quotients $\mathbb{Z}_2\wr\mathbb{Z}_{n_k}$ furnish the chain, where $(n_k)_{k\in \N}$ is any computable sequence tending to infinity and satisfying $n_k\mid n_{k+1}$ for all $k\in \N$.
\end{enumerate}
\end{examples}

We note that there are examples of computable presentations of amenable groups for which Corollary \ref{cor:rf} does not necessarily apply and yet the corresponding presentation of the group \cstar-algebra is effectively Kirchberg nuclear.  For example, if $G^\#$ is a computable presentation of an abelian group, then $C^*(G)^\#$ is computable by \cite[Corollary 3.22]{fox2022computable}; since $C^*(G)$ is abelian, $C^*(G)^\#$ is effectively Kirchberg nuclear by Corollary \ref{cor:abeliankirch}.  Similarly, if $G^\#$ is a computable presentation of a locally finite group, then $C^*(G)^\#$ is easily seen to be a computable presentation; since $C^*(G)$ is AF, we obtain that $C^*(G)^\#$ is effectively Kirchberg nuclear by Proposition \ref{prop:AFkirch}.

\section{$\Pi_2^0$-completeness}\label{sec:completeness}

\newcommand{\amen}{\operatorname{Amen}}
\newcommand{\Inf}{\textbf{Inf}}
\newcommand{\group}{\operatorname{Grp}}
\newcommand{\Nuc}{\operatorname{Nuc}}
\newcommand{\set}[2]{\ensuremath{ \{ #1 : #2 \} }}

Set $\Nuc$ to be the set of indices of computable presentations of nuclear \cstar-algebras.  We begin with the following:

\begin{proposition}\label{prop:nucispi2}
$\Nuc$ is $\Pi_2^0$.
\end{proposition}

\begin{proof}
 The set of indices of computable presentations of \cstar-algebras is itself a $\Pi_2^0$-set (see \cite[Theorem 12.1]{AFPaper}).  Moreover, the proof of the direction (1) implies (2) in Theorem \ref{mainequivalence} shows that being nuclear is indeed a $\Pi_2^0$-condition on computable presentations of \cstar-algebras.  Indeed, letting $X_k$ denote the set of tuples $(\vec a,n,\vec b,w)$ satisfying the three bullet points appearing in Theorem \ref{mainequivalence} (where $\vec a$ is a $k$-tuple), we see that $X_k$ is computable uniformly in the presentation $A^\#$; the proof of Theorem \ref{mainequivalence} shows that the algebra is nuclear if and only if, for each $k$ and each $k$-tuple $\vec a$, there is a triple $(n,\vec b,w)$ such that $(\vec a,n,\vec b,w)\in X_k$, establishing that being nuclear is $\Pi_2^0$.
$\Nuc$ is thus the intersection of 
$\Pi^0_2$-sets, hence $\Pi^0_2$.
\end{proof}

The motivating question considered in this section is the following:

\begin{question}\label{completequestion}
Is $\Nuc$ $\Pi^0_2$-complete within the class of indices of computable presentations of \cstar-algebras?
\end{question}

Here $\Pi^0_2$-completeness \emph{within the class} means that for every $\Pi^0_2$ set $A$, there is a computable total function $f$ mapping $\omega$ \emph{into the class of indices of computable presentations of \cstar-algebras} such that
$n\in A\iff f(n)\in\Nuc$.  Contrast this with the requirement for ordinary $\Pi^0_2$-completeness, 
where one would only need $f:\omega\to\omega$.  
Ordinary $\Pi^0_2$-completeness of $\Nuc$ (which is not difficult to prove)
would imply only that being the index of a nuclear \cstar-algebra -- as opposed to being the index of something that is maybe not even a \cstar-algebra at all -- is $\Pi^0_2$-hard.
$\Pi^0_2$-completeness within the class of indices of computable presentations of \cstar-algebras, if true, would mean that even as a property of \cstar-algebras, nuclearity is still $\Pi^0_2$-hard.

There is some precedence for a positive answer to Question \ref{completequestion}.  In \cite[Theorem 12.4]{AFPaper}, it was proven that the set of indices for computable presentations of UHF algebras is $\Pi_2^0$-complete, whereas in \cite[Theorems 12.2 and 12.5]{AFPaper}, it was shown that the set of indices for computable presentations of (simple) AF-algebras is $\Pi_2^0$-complete.

In this section, we describe partial progress towards settling Question \ref{completequestion}.  We begin with a detour giving a positive answer to the analogous question for groups, which we believe is of independent interest.

\subsection{Amenability}

First consider the set 
$$\group=\set{e\in\omega}{\phi_e\text{~computes a group on the domain~}\omega}.$$
Here we consider $\phi_e$ as a function
 $\phi_e:\omega^2\to \omega$ and require this binary function to be total and to satisfy the group axioms.

\begin{lemma}
$\group$ is a $\Pi_2^0$ set.
\end{lemma}

\begin{proof}
It suffices to note that $\group$ is the intersection of $\operatorname{Tot}:=\{e\in \omega \ : \ \phi_e \text{ is total}\}$ with the following sets:
\begin{itemize}
    \item $\{e\in \omega \ : \ (\exists x\in \omega)\phi_e(x,x)=x\}$.
    \item $\{e\in \omega \ : \ (\forall x,y\in \omega)((\phi_e(x,x)=x\wedge \phi_e(y,y)=y)\rightarrow x=y)\}$.
    \item $\{e\in \omega \ : \ (\forall x\in \omega)(\exists y\in \omega)(\phi_e(\phi_e(x,y),\phi_e(x,y))=\phi_e(x,y)\\ \wedge \phi_e(\phi_e(y,x),\phi_e(y,x))=\phi_e(y,x)\}$.
    \item $\{e\in \omega \ : \ (\forall x,y,z\in \omega)\phi_e(\phi_e(x,y),z)=\phi_e(x,\phi_e(y,z))\}$.\qedhere
\end{itemize}
\end{proof}

Now consider the set
$$\amen=\set{e\in\omega}{\phi_e\text{~computes an amenable group on the domain~}\omega}.$$

\begin{lemma}
$\amen$ is a $\Pi_2^0$ set.
\end{lemma}

\begin{proof}
Recall that amenability can be characterized by the following condition:
    $$ (\forall\text{~finite~}F_0\subseteq G)(\forall\epsilon>0)(\exists\text{~finite~}F_1\subseteq G)
(\forall x\in F_0)\frac{|xF_1\bigtriangleup F_1|}{|F_1|}<\epsilon.$$
By standard coding tricks, this is easily seen to be a $\Pi_2^0$-condition.
\end{proof}

We now prove the main result of this subsection.

\begin{theorem}
\label{thm:amenability}
$\amen$ is a $\Pi_2^0$-complete set.
Indeed, it is $\Pi^0_2$-complete
within the class $\group$.
\end{theorem}

\begin{proof}[Proof of Theorem \ref{thm:amenability}]
We define a computable total function $R:\omega\to \group$ such that 
$$(\forall e)~[e\in \Inf \iff R(e)\in\amen].$$
Here $\Inf:=\set{e}{W_e\text{~is infinite}}$.
$R$ can then be converted into an injective total computable function, by the Padding Lemma,
showing the $\Pi^0_2$-complete set $\Inf$ to be $1$-reducible to $\amen$ and thus giving the result.  Further details
on all of these concepts may be found in \cite[chs.\ I-IV]{Soare}.

  In what follows, a \emph{fragment} of a group $G$ is a substructure of $G$ in the relational signature with a ternary predicate defining the group operation.  (Thus, a fragment retains the group's structure but need not be closed under the group operation.)

For a given $e\in\omega$, we define $R(e)$ to be the index of the group operation
$\cdot$ on the domain $\omega$ in the computable group $F$ that we
now build effectively as the increasing union of group fragments $F_s$ at each stage $s$.
Every time a new element appears in $W_e$, say at stage $s$, we will see this and will make
$F_{s+1}$ into a finite group.  At stages when no new element appears in $W_e$, we will
take one step towards building a copy of $FG_2$, the free group on two generators, within $F$.
If $e\notin\Inf$, then $F$ will contain a copy of $FG_2$, built without interruption starting
at the stage after the last element of $W_e$ appeared, and hence will
not be amenable (so $R(e)\notin\amen$).  On the other hand, if $e\in\Inf$, then
$F$ will be a nested union of finite groups $F_{s_1}\subseteq F_{s_2}\subseteq\cdots$,
and therefore amenable.

We start with $F_0$ containing the first three domain elements $0,1,2$ from $\omega$.
We define $U_0$ to be the trivial group and fix a computable
presentation $G_0=U_0\times FG_2$ with two specific generators $(0,a_0)$
and $(0,b_0)$ of the $\mathbb{F}_2$ factor,  (The first factor may become nontrivial
as the construction proceeds.)  Then we
create the fragmentary embedding $f_0:F_0\hookrightarrow G_0$ with $f_0(0)=(0,\lambda)$
the identity in $G_0$, $f_0(1)=(0,a_0)$, and $f_0(2)=(0,b_0)$.  This embedding
will serve as a guideline to our construction of $F$ as long as $W_e$ stays empty.

At each stage $s+1$, we know the fragment $(F_s,\cdot)$ and the current embedding $f_s$,
which maps $F_s$ into a group $G_s\cong U_s\times FG_2$, where $U_s$ is finite.
If $W_{e,s+1}=W_{e,s}$, then we keep $U_{s+1}=U_s$ and $G_{s+1}= G_s$,
find the least element $y_{s+1}\in G_{s+1}$ that is not yet in the (finite) image of $f_s$, define $f_{s+1}=f_s\cup (x_{s+1},y_{s+1})$, where $x_{s+1}=|F_s|$ is the least
number in $\omega-F_s$, and extend $\cdot$ to $F_{s+1}$ so that this $f_{s+1}$
remains a fragmentary embedding.  Thus we have taken one more step towards building
a copy of $FG_2$ within $F$.

If $W_{e,s+1}\neq W_{e,s}$, we find a finite group $V_{s+1}$
such that the finite image $f_s(F_s)$, when projected onto the second coordinate
$\mathbb{F}_2$ of $G_s$, can be extended to a copy of $V_{s+1}$ via some group
homomorphism $g_{s+1}:\mathbb{F}_2\to V_{s+1}$ that is injective
on the subdomain $\pi_1(f_s(F_s))$.  The heart of this proof is Lemma \ref{lemma:infstep}
and its corollary below, which shows that such a $V_{s+1}$ and such a $g_{s+1}$ exist.
(Since they exist, we can find such $V$ and $g_{s+1}$ just by searching systematically:
when we stumble across them, we will recognize that they suffice.)
Using $g_{s+1}$, we extend $F_s$ to a copy $F_{s+1}$ of $U_s\times V_{s+1}$
with an isomorphism $f_{s+1}:F_{s+1}\to U_s\times V_{s+1}$ such that for each $x\in F_s$,
$$f_{s+1}(x) = ((\pi_0(f_s(x)), g_{s+1}(\pi_1(f_s(x)))),\text{id})\in (U_s\times V_{s+1})\times FG_2.$$
(Here $\pi_0$ and $\pi_1$ represent projection onto the first and second coordinates, respectively.)
We define $U_{s+1}=U_s\times V_{s+1}$ and create the new computable group
$G_{s+1}=U_{s+1}\times \mathbb{F}_2$.
Notice that this step does not redefine the (partial) group operation $\cdot$ on $F_s$
in any way:  we are continuing to construct a computable group, with
no adjustment of the portion already constructed..  This completes stage $s+1$.

Since the notation above is dense, we summarize here the action when $W_{e,s+1}\neq W_{e,s}$.
First we find a finite group $V_{s+1}$ such that the portion of $FG_2$ already
used (as second coordinates) in the finite subset $F_s$ of $G_s$ can be
consistently extended to a copy of $V$.  Having found it (by blind search,
with an appeal to Corollary \ref{cor:infstep}), we extend the currently-defined
finite portion $F_s$ of $U_s\times FG_2$ consistently to a full copy of $U_s\times V$,
adding fresh elements to it to make $F_{s+1}\cong U_s\times V$.
This $U_s\times V$ becomes the new $U_{s+1}$, which is finite, and the new ``target group''
$G_{s+1}$ is the direct product of this $U_{s+1}$ with an entirely new copy of $FG_2$.
The new embedding $f_{s+1}$ maps $F_{s+1}$ onto $U_{s+1}\times\{\lambda\}$ within $G_{s+1}$.

Now we show that $R$ is an $m$-reduction.
If $e\in\Inf$, then there are infinitely many stages $s_1<s_2<\cdots$ at which
$W_{e,s_i}\neq W_{e,s_i-1}$.  (By convention each $W_{e,s}\subseteq\{ 0,1,\ldots,s\}$,
so only finitely many elements enter $W_e$ at any single stage.)
Consequently, for each stage in this sequence, $F_{s_i}\cong F_{s_{i-1}}\times V_{s_i}$,
with $F_{s_{i-1}}$ embedding into $F_{s_i}$ as the first factor. 
The group $F$ that we construct across all $\omega$-many stages
is thus isomorphic to the algebraic direct product
$V_{s_1}\oplus V_{s_2}\oplus\cdots$.  It follows that $F$ is locally amenable,
as every finite subset of $F$ lies within some finite
(hence amenable) group $V_{s_1}\oplus\cdots\oplus V_{s_k}$,
and consequently $F$ is amenable.

Conversely, if $e\notin\Inf$, let $t$ be minimal (possibly $0$) such that $W_{e,t}=W_e$.
Then at every stage $s+1>t$, we simply extended $F_{s+1}$ by one more
element, defining the group operation so that $f_{s+1}$ can map
the new element to the least available element of $G_{s_k}$
and remain a partial group embedding.  In the limit, $f=\lim_s f_s$
will thus be a group isomorphism from $F=\bigcup_s F_s$
onto $G=G_t=U_{t}\times FG_2$.  Thus $F$ contains
a copy of $FG_2$ and consequently cannot be amenable.
This completes our proof, modulo Corollary \ref{cor:infstep}.
\end{proof}

The key lemma for this theorem was proven by Iwasawa in \cite[Satz 2]{Iwasawa}
(in German, translated here).
\begin{lemma}[Iwasawa \cite{Iwasawa}]
\label{lemma:infstep}
Let $F$ be a free group on any number of generators.  Then the intersection
of all finite-index normal subgroups of $F$ is the trivial subgroup of $F$.
\end{lemma}
\begin{corollary}
\label{cor:infstep}
Whenever $F_0$ is a finite fragment of the group $FG_2$  (or indeed of any free group),
there exists a finite group $V$ such that the (finite, partial) atomic diagram of $F_0$
within $F$ extends to the atomic diagram of a copy of $V$.
\end{corollary}
\begin{proof}
To prove this, we need a finite-index normal subgroup $N\triangleleft FG_2$
such that $N\cap (F_0\cdot F_0^{-1})=\{\lambda\}$, where $\lambda$ is the empty word
in $\mathbb{F}_2$.  Then $G/N$
will be finite and the quotient map will send each $x\neq y$ in $F_0$ to distinct
cosets (since $xy^{-1}\notin N$), so the atomic diagram of $G/N$ will realize
the Corollary.

By Lemma \ref{lemma:infstep}, for each non-identity $z\in F_0F_0^{-1}$,
there is a finite-index normal subgroup $N_z\triangleleft FG_2$ with $z\notin F_z$.
Their intersection
$$N=\bigcap_{z\in F_0F_0^{-1},z\neq\lambda} N_z$$
is then also normal of finite index in $FG_2$ and contains none of the proscribed elements.
\end{proof}

\subsection{Inner amenability}

Although it will not be relevant to our study of nuclear \cstar-algebras, we add here a strengthening of Theorem \ref{thm:amenability}.  Set $$\operatorname{InnAmen}:=\{e\in \omega \ : \ \phi_e \text{ computes an inner amenable group on the domain }\omega\}.$$
This is a superset of $\amen$, as all amenable groups are inner amenable.

\begin{lemma}
$\operatorname{InnAmen}$ is a $\Pi_2^0$ set.
\end{lemma}

\begin{proof}
It suffices to recall that inner amenability also admits a Folner-type criteria, namely $G$ is inner amenable if and only if
$$ (\forall\text{~finite~}F_0\subseteq G)(\forall\epsilon>0)(\exists\text{~finite~}F_1\subseteq G)
(\forall x\in F_0)\frac{|xF_1x^{-1}\bigtriangleup F_1|}{|F_1|}<\epsilon.\qedhere$$
\end{proof}

\begin{theorem}
\label{prop:inneramenable}
$\operatorname{InnAmen}$ is $\Pi^0_2$-complete within the class $\group$.
\end{theorem}
\begin{proof}
The construction is quite similar to that in the proof of Theorem \ref{thm:amenability}.  The only difference is that here, instead of trying (at stages $s+1$ where $W_{e,s+1}=W_{e,s}$) to build a copy of
the free group $FG_2$ within our computable group $G$,
we will try instead to build $G$ isomorphic to the group $\operatorname{SL}_2\left(\overline{\F_2(t)}\right)$, which is known (see \cite[\S 2]{Goldbring.2024}) not to be inner amenable.
(Here $\F_2$ is the two-element field, which we extend by a transcendental element $t$ and then take the algebraic closure $\overline{\F_2(t)}$.)
Now if any stage $s$ exists for which $W_e=W_{e,s}$ (that is, if $W_e$ is finite),
then from stage $s+1$ on, our construction will succeed in making
$G$ isomorphic to $SL_2\left(\overline{\F_2(t)}\right)$.
On the other hand, if $W_e$ is infinite, the construction will produce a $G$ that is a direct limit of finite subgroups of itself, hence amenable and thus inner amenable.

The construction begins with $K_0=\F_2$ and $G_0=\operatorname{SL}_2(K_0)$,
with the latter viewed as a finite fragment of (a fixed computable presentation of) the matrix group $\operatorname{SL}_2\left(\overline{\F_2(t_0)}\right)$, where $t_0$ is transcendental over $\F_2$.
We now explain the procedure at each stage $s+1$.  At such a stage, we begin with a finite subfield $K_s$ of the algebraic closure $\overline{\F_2}$ and a finite fragment $G_s$ of the group $\operatorname{SL}_2\left(\overline{K_s(t_s)}\right)$.  If $W_{e,s+1}=W_{e,s}$, then we simply extend $G_{s+1}$ to include the next matrix in $\operatorname{SL}_2\left(\overline{\F_2(t_s)}\right)$, setting $t_{s+1}=t_s$ and $K_{s+1}$ equal to the finite subfield of $\overline{\F_2}$ generated by all elements of $\overline{\F_2(t_{s+1})}$ that appear in the matrices now in $G_{s+1}$ and are algebraic over $\F_2$.  Thus, if in fact $W_e=W_{e,s+1}$, then the execution of this step at all subsequent stages will keep $t_{s+1}$ transcendental over $\F_2$ forever and build $G=\bigcup_n G_n=\operatorname{SL}_2\left(\overline{\F_2(t_{s+1})}\right)$ as desired.

If $W_{e,s+1}\neq W_{e,s}$, however, then we wish to increase the current fragment $G_s$ of $\operatorname{SL}_2\left(\overline{K_s(t_s)}\right)$
to a finite group.  $G_s$ uses only finitely many elements of the algebraic closure $\overline{\F_2(t_s)}$, so we can effectively convert the current $t_s$ into a field element algebraic over $\F_2$,
in a way consistent with the current field fragment.  (That is, we embed the current field fragment into the algebraic closure of $\F_2$.)  This is possible because the type of a transcendental over $\F_2$ is nonprincipal, hence cannot be generated by the finitely many atomic facts given by the current field fragment.
Since the current fragment is finite, so is its image in $\overline{\F_2}$,
and therefore so is our new $K_{s+1}$, which is the subfield generated in $\overline{\F_2}$ by the current fragment and the new algebraic element $t_s$.  Next we use this field-fragment embedding to embed the current finite group fragment $G_s$ of 
$\operatorname{SL}_2\left(\overline{\F_2(t)}\right)$
into $\operatorname{SL}_2(K_{s+1})$ and then extend it to the full (finite) group $\operatorname{SL}_2(K_{s+1})$, defining $G_{s+1}=\operatorname{SL}_2(K_{s+1})$.  Thus $G_s$, which was a fragment of $\operatorname{SL}_2\left(\overline{\F_2(t_s)}\right)$,
has now been enlarged to an actual finite group $G_{s+1}$, not just a fragment.
Therefore, if $W_e$ increases infinitely often, then $G$ will turn out to be the limit of the sequence
$G_0\hookrightarrow G_{s_1}\hookrightarrow G_{s_2}\hookrightarrow\cdots$ of finite groups, as desired, where
$s_1<s_2<\cdots$ are the stages at which $W_e$ acquires new elements.
\end{proof}

Since amenability implies inner amenability,
Theorem \ref{thm:amenability} actually follows directly from the proof of Theorem \ref{prop:inneramenable}.  We included our proof of Theorem \ref{thm:amenability} on grounds that it is both simpler and purely group-theoretic, while the latter proof required substantial use of basic field theory as well.

\subsection{Back to Question \ref{completequestion}}

Let $\Nuc'$ denote the set of indices of \emph{computably enumerable} presentations of nuclear \cstar-algebras.

\begin{corollary}\label{cor:nuchard}
$\Nuc'$ is $\Pi_2^0$-hard within the class of indices of c.e.\ presentations of \cstar-algebras.
\end{corollary}

\begin{proof}
    Let $h:\omega\to \omega$ be the following computable function:  given $e\in \omega$, first compute $R
(e)$, where $R:\omega\to \group$ is the computable injective function defined in the proof of Theorem \ref{thm:amenability}.  Let $G_e$ denote the group on $\omega$ whose group operation is computed by $\phi_{R(e)}$.  Now let $h(e)$ be the index for the ``induced'' c.e.\ presentation of the universal group \cstar-algebra $C^*(G_e)$.  More precisely, $C^*(G_e)$ is the universal \cstar-algebra on countably many generators subject to the relations given by the multiplication table $\phi_{R(e)}$ and such presentations admit canonical c.e.\ presentations as described in \cite[Theorem 3.3]{fox2022computable}.  It remains to show that, for all $e\in \omega$, we have that $e\in \Inf$ just if $h(e)\in \operatorname{Nuc}'$.  On the one hand, if $e\in \Inf$, then $R(e)\in \operatorname{Amen}$, whence $C^*(G_e)$ is nuclear.  Moreover, by \cite[Corollary 3.22]{fox2022computable}, the presentation coded by $h(e)$ is in fact computable, whence $h(e)\in \operatorname{Nuc}'$.  On the other hand, if $e\notin \Inf$, then $G_e$ is not amenable, whence $C^*(G_e)$ is not nuclear and thus $h(e)\notin \operatorname{Nuc}'$, yet
$h(e)$ is still the index of a computable (non-nuclear) \cstar-algebra.
\end{proof}

\begin{question}\label{ques:nuc'pi2}
Is $\Nuc'$ $\Pi_2^0$?
\end{question}

In light of Corollary \ref{cor:nuchard}, a positive answer to this question would show that $\Nuc'$ is $\Pi_2^0$-complete within the class of indices of c.e.\ presentations of \cstar-algebras, lending some credence to a positive answer to Question \ref{completequestion}.

In \cite[Theorem 12.1]{AFPaper}, it was remarked that the set of indices of c.e.\ presentations of \cstar-algebras is $\Pi_2^0$.  On the other hand, note that the proof that $\Nuc$ is $\Pi_2^0$ used computability of the presentation in an essential way via its use of the effective Arveson extension.  
\bibliographystyle{amsplain}
\bibliography{bib}

\end{document}